%Paper: alg-geom/9509005
%From: "Steven L. Kleiman" <kleiman@math.mit.edu>
%Date: Sun, 10 Sep 1995 00:06:32 -0400

%%
%% AmS-TeX- Version 2.1
%% AMS-TeX- documentstyle: amsppt.sty version="2.1c" 11-Dec-1992
%%              also compiles with version "2.1a"
%
%% SOME CHOICES
\input amstex.tex              % Maybe eliminated if part of format file
\documentstyle{amsppt}         % Maybe eliminated if part of format file
\nologo                         % Eliminates \AMSTeX logo from first page
\def\TheMagstep{\magstep1}	% Normal magnification
		% Changed to \magstep0 by \DoublepageOutput{TRUE}
\def\PaperSize{letter}		% \PaperSize is used to
% \def\PaperSize{AFour}		%  center text on page
%% To get two pages side by side in landscape mode on an ordinary
%%    sheet via dvips and  PostScript, enable the next command.
% \def\DoublepageOutput{TRUE}
%%	END of CHOICES

 \advance\vsize-3\normalbaselineskip
\def\papernumber{alg-geom/9509005}
\def\thedate{950911}
%\input sym.mac			%% Macro file
%%  MATH MACROS  %%
\let\:=\colon
\let\To=\longrightarrow 
\def\onto{\to\mathrel{\mkern-15mu}\to}
\def\risom{\,\overset{\sim}\to{\smash\longrightarrow\vrule
    height0.5ex width0pt depth0pt}\,}
\let\?=\overline  
\let\impls=\Rightarrow
\let\ox=\otimes\let\x=\times
\let\m=m	% \def\m{{\bold m}}
\let\D=\Delta \def\O{\Cal O}
\def\Ser#1{{\(S${}_{#1}$)}}
\def\ext#1{\hbox{\rm Ext}^{#1}_} \def\hom{\hbox{\rm Hom}}
\def\I#1{I^{(#1)}} \def\BI{\bold{I}}
\def\and{\hbox{ \rm and }} 
\def\SR{\mathop {{\Cal R}_s}}
\def\up#1{\raise 3pt \hbox{$#1$}}
\def\col #1 #2 #3 {\pmatrix #1\\ #2\\ #3\\\endpmatrix}
\let\wt=\widetilde

\def\mathopdef#1{\expandafter\def
 \csname#1\endcsname{\operatorname {#1}}}
\def\NoOp{*!*}
\def\NextOp#1 {\def\TheOp{#1}\ifx\TheOp\NoOp\let\next\relax
  \else\mathopdef{#1}\let\next\NextOp \fi \next}
\NextOp
 Hom adj grade Fitt Spec Ann Ass depth Ext height gr Cok Im Ker char
rank length cod
 *!*
\let\rk=\rank

%% Display numbers
\def\dno#1${\eqno\hbox{\rm(\number\sectno.#1)}$}
\def\Cs#1){\unskip~{\rm(\number\sectno.#1)}}
%% Convient macro for display numbers
\newcount\displayno
\def\eqlt#1${\global\advance\displayno by 1
 \expandafter\xdef
  \csname \the\sectno#1\endcsname{\the\displayno}
 \eqno(\the\sectno.\the\displayno)$}
\def\disp#1#2{(#1.\csname#1#2\endcsname)}
\def\Cn#1){\unskip~{\rm(\number\sectno.\csname\the\sectno#1\endcsname)}}

%% FONTS
% \font\twelvebf=cmbx12
      %% STYLE %%

%% For Roman paranthetical material in nonRoman text
\def\(#1){{\rm(#1)}}\let\leftp=(
\def\activeleftp{\catcode`\(=\active}
{\activeleftp\gdef({\ifmmode\let\next=\leftp \else\let\next=\(\fi\next}}
\def\biglp{\bigl(}

%% SECTIONING

%% Articles
\def\artkey #1 {{\bf (\number\sectno.#1)}\enspace}

\def\eg  {\medbreak {\bf Example }\artkey}
\def\proclaim#1 #2 {\medbreak{\bf#1 }\artkey#2 \bgroup\it\activeleftp}
\def\endproclaim{\egroup\medbreak}
\def\pf{\endproclaim{\bf Proof.}\enspace}
\def\prop{\proclaim Proposition } \let\prp=\prop
\def\lem{\proclaim Lemma }
\def\thm{\proclaim Theorem }
\def\cor{\proclaim Corollary }

%% To start a new section.
%%  Checks that there's room enough on the page.
 \newcount\sectno \sectno=0
 \newskip\sectskipamount \sectskipamount=0pt plus30pt
 \def\newsect #1\par{\displayno=0 %\proclaimno=0
   \advance\sectno by 1
   \vskip\sectskipamount\penalty-250\vskip-\sectskipamount
   \bigskip%\smallskip \medskip
   \centerline{\bf \number\sectno. #1}\nobreak
   \medskip%\smallskip
   \message{#1 }
}

%% Display numbers
\def\dno#1${\eqno\hbox{\rm(\number\sectno.#1)}$}
\def\Cs#1){\unskip~{\rm(\number\sectno.#1)}}
%% Convient macro for display numbers
\newcount\displayno
\def\eqlt#1${\global\advance\displayno by 1
 \expandafter\xdef
  \csname \the\sectno#1\endcsname{\the\displayno}
 \eqno(\the\sectno.\the\displayno)$}
\def\disp#1#2{(#1.\csname#1#2\endcsname)}
\def\Cn#1){\unskip~{\rm(\number\sectno.\csname\the\sectno#1\endcsname)}}

       %% FORMAT
%% Page layout
\def\TRUE{TRUE}
\ifx\DoublepageOutput\TRUE \def\TheMagstep{\magstep0} \fi
\mag=\TheMagstep
\TagsOnRight
\parskip=0pt plus3pt
\abovedisplayskip 6pt plus3pt minus3pt
\belowdisplayskip=\abovedisplayskip
\nopagenumbers

% Center text on page
	% additional vertical adjustment
\newskip\vadjustskip \vadjustskip=0.5\normalbaselineskip
\def\centertext
 {\hoffset=\pgwidth \advance\hoffset-\hsize
  \advance\hoffset-2truein \divide\hoffset by 2
  \voffset=\pgheight \advance\voffset-\vsize
  \advance\voffset-2truein \divide\voffset by 2
  \advance\voffset\vadjustskip
 }
\newdimen\pgwidth\newdimen\pgheight
\def\letter{letter}\def\AFour{AFour}
\ifx\PaperSize\letter
 \pgwidth=8.5truein \pgheight=11truein
 \message{- Got a paper size of letter.  }\centertext
 \fi
\ifx\PaperSize\AFour
 \pgwidth=210truemm \pgheight=297truemm
 \message{- Got a paper size of AFour.  }\centertext \fi

%% Two-column landscape format
% Modified from the TeX book, p. 257.
 \newdimen\fullhsize \newbox\leftcolumn
 \def\fulline{\hbox to \fullhsize}
\def\doublepageoutput
{\let\lr=L
 \output={\if L\lr
          \global\setbox\leftcolumn=\columnbox \global\let\lr=R%
        \else \doubleformat \global\let\lr=L\fi
        \ifnum\outputpenalty>-20000 \else\dosupereject\fi}%
 \def\doubleformat{\shipout\vbox{%
        \fulline{\hfil\hfil\box\leftcolumn\hfil\columnbox\hfil\hfil}%
}%
 }%
 \def\columnbox{\vbox
   {\makeheadline\pagebody\makefootline\advancepageno}%
   }
 \fullhsize=\pgheight \hoffset=-1truein
 \voffset=\pgwidth \advance\voffset-\vsize
  \advance\voffset-2truein \divide\voffset by 2
  \advance\voffset\vadjustskip
\headline={%
 \eightpoint
  \ifnum\pageno=1\firstheadline
  \else
    \ifodd\pageno\rightheadline
    \else\leftheadline
    \fi
  \fi
}
\let\firstheadline=\hfil

%\null\vfill\nopagenumbers\eject\pageno=1\relax % to put page on right
}
\ifx\DoublepageOutput\TRUE \doublepageoutput \fi

%THE DATE:
\def\today{\number\day \space\ifcase\month\or
 January\or February\or March\or April\or May\or June\or
 July\or August\or September\or October\or November\or December\fi
\space \number\year}

%% HEADLINES %%
 \leftheadtext{KLEIMAN AND ULRICH}
 \rightheadtext{Gorenstein rings, symmetric matrices, etc.}
 \def\rightheadline{\rlap{\thedate}\hfil
\the\rightheadtoks\hfil\llap{\folio}}
 \def\leftheadline{\rlap{\folio}\hfil
 \the\leftheadtoks\hfil\llap{\papernumber}}

%% Redefine \item to give the indentation of AMSTeX and the roman font
\def\item#1 {\par\indent\indent\indent\indent
 \hangindent4\parindent
 \llap{\rm (#1)\enspace}\ignorespaces}
%% Define a similar macro without the hanging indentation for assertions
%% and  that starts each part with an ordinary \parindent
 \def\part#1 {\par{\rm (#1)\enspace}\ignorespaces}

%% REFERENCING %%
%%  Macros for introducing the reference keys in order
 \newif\ifproofing \proofingfalse % Decide between Alpha and Num keys.
 \newcount\refno \refno=0
 \def\MakeKey{\advance\refno by 1 \expandafter\xdef
 \csname\TheKey\endcsname{{%
\ifproofing\TheKey\else\number\refno\fi}}\NextKey}
 \def\NextKey#1 {\def\TheKey{#1}\ifx\TheKey\NoKey\let\next\relax
  \else\let\next\MakeKey \fi \next}
 \def\NoKey{*!*}
 \def\RefKeys #1\endRefKeys{\expandafter\NextKey #1 *!* }
\def\SetKey#1{{\bf\csname#1\endcsname}}

%% \cite code from AMSTeX, modified so first argument is a control sequence
\catcode`\@=11
\def\relaxnext@{\let\next\relax}
\def\cite#1{\relaxnext@
 \def\nextiii@##1,##2\end@{{\rm[\SetKey{##1}, \sfcode`\.=1000##2]}}%
 \in@,{#1}\ifin@\def\next{\nextiii@#1\end@}\else
 \def\next{\unskip\space{\rm[\SetKey{#1}]}}\fi\next}
\newif\ifin@
\def\in@#1#2{\def\in@@##1#1##2##3\in@@
 {\ifx\in@##2\in@false\else\in@true\fi}%
 \in@@#2#1\in@\in@@}

%%  ``Ties'' with a \thinspace for page numbers
\def\p.{\unskip\space p.\penalty\@M \thinspace}
\def\pp.{\unskip\space pp.\penalty\@M \thinspace}
%% Roman parentheses for use within text set in other fonts
%% taken from the TeXbook, p.409, but modified to use the italic correction
%\def\(#1){{\rm(}#1\/{\rm)}}

\catcode`\@=13

%% Macros for setting references
 \def\SetRef#1 #2\par{\rm
   \hang\llap{[\csname#1\endcsname]\enspace}%
   \ignorespaces#2\unskip.\endgraf}
 \newbox\keybox \setbox\keybox=\hbox{[18]\enspace}
 \newdimen\keyindent \keyindent=\wd\keybox
 \def\references{\vskip-\smallskipamount
  \bgroup   \eightpoint   \frenchspacing
   \parindent=\keyindent  \parskip=\smallskipamount
   \everypar={\SetRef}}
 \def\endreferences{\egroup}
 \def\nocomma{\def\eatcomma##1,{}\expandafter\eatcomma}

%% Papers and books
  \def\paper{\unskip, \bgroup\it}
 \def\paperinfo{\unskip, \rm}
  \def\inbook#1\bookinfo#2\publ#3\yr#4\pages#5
   {\unskip, \egroup in ``#1\unskip,'' #2\unskip, #3\unskip,
   #4\unskip, pp.~#5}
 \def\at.{.\spacefactor3000}
 \def\book{\unskip, \rm ``}
 \def\bookinfo#1{\unskip," #1}
 \def\preprint{\unskip, \egroup\rm preprint, }
 \def\mams #1 #2 {\unskip, \egroup Mem. Amer. Math. Soc. {\bf #1} (#2)}
%% SERIALS
 \def\serial#1#2{\expandafter\def\csname#1\endcsname ##1 ##2 ##3
  {\unskip, \egroup #2 {\bf##1} (##2), ##3}}
 \serial{ajm}{Amer. J. Math.}
 \serial{cras}{C. R. Acad. Sci. Paris}
 \serial{cmh}{Comment. Math. Helv.}
 \serial{comp}{Comp. Math.}
 \serial{ja}{J. Algebra}
 \serial{jmsj}{J. Math. Soc. Japan}
 \serial{jmku}{J. Math. Kyoto Univ.}
 \serial{jram}{J. reine angew. Math.}
 \serial{lmslns}{London Math. Soc. Lecture Note Series}
 \serial{ma}{Math. Ann.}
 \serial{mpcps}{Math. Proc. Camb. Phil. Soc.}
 \serial{nmj}{Nagoya Math J.}
 \serial{pams}{Proc. Amer. Math. Soc.}
 \serial{prs}{Proc. Royal Soc.}
 \serial{tams}{Trans. Amer. Math. Soc.}

 \def\patzcuaro{\inbook Algebraic geometry and complex analysis
\bookinfo E. Ram\'\i rez de Arellano (ed.), Proc. Conf., P\'atzcuaro
1987, Lecture Notes in Math. {\bf 1414} \publ Springer-Verlag \yr1989
\pages}
 \def\arcireale{\inbook Complete Intersections
\bookinfo S. Greco and R. Strano (eds.), Proc. Conf., Arcireale 1983
Lecture Notes in Math. {\bf 1092} \publ Springer-Verlag \yr1984
\pages}
 
 \def\bucharest{\inbook Algebraic Geometry, Bucharest 1982\at.
\bookinfo L. B\v adescu and D. Popescu
(eds.), Lecture Notes in Math. {\bf 1056} \publ Springer-Verlag
\yr 1984 \pages}

\RefKeys
 AN Brod Bruns BE77a BE77b Cat dJvS EN Eisenbud Ei-Ma Flenner G-N G-N-S
G-N-W Grassi Groth HeK H-U Jzfk JP KLU92 KLU Kunz Kutz Mats MP Reiten
Roberts85 Roberts90 Shamash VV Valla84
 \endRefKeys

%%% DOCUMENT
 \topmatter
 \title
	Gorenstein algebras, symmetric matrices,
self-linked ideals, and symbolic powers
 \endtitle
 \author
	Steven KLEIMAN,$^1$ and
	Bernd ULRICH$^2$
 \endauthor
 \address
   Mathematics Department Room 2-278, Massachusetts Institute of
Technology,
   Cam\-bridge, MA 02139-4307, USA.
 \endaddress
 \email \tt Kleiman\@math.MIT.edu \endemail
 \address
  Mathematics Department, Michigan State University, East Lansing,
MI 48824-1027, USA.
 \endaddress
 \email \tt  Ulrich\@math.MSU.edu \endemail

\date September 11, 1995\enddate
 \thanks
 {\it Acknowledgements.}\enspace It is a pleasure to thank David Mond
for many stimulating discussions.  It is a pleasure to thank David
Eisenbud for his comments, which led to the present improved versions
of (2.1) and (2.2).  It is a pleasure to thank Hubert Flenner for his
help with Bertini's theorem for local rings.  Finally, it is a pleasure
to thank Joseph Lipman for his long collaboration on the present
subject; in the estimation of the two authors de jure, if not in his
own, he is a third author de facto.  \endgraf
 $^1$Supported in part by NSF grant 9400918-DMS.  It is a pleasure for
this author to thank the Mathematical Institute of the University of
Copenhagen for its hospitality during the summer of 1995 when this work
was completed.\endgraf
	    $^2$Supported in part by NSF grant DMS-9305832.
 \endthanks
% \keywords  \endkeywords
 \abstract
	Inspired by recent work in the theory of central projections
onto hypersurfaces, we characterize self-linked perfect ideals of grade
2 as those with a Hilbert--Burch matrix that has a maximal symmetric
subblock.  We also prove that every Gorenstein perfect algebra of grade
1 can be presented, as a module, by a symmetric matrix.  Both results
are derived from the same elementary lemma about symmetrizing a matrix
that has, modulo a nonzerodivisor, a symmetric syzygy matrix.  In
addition, we establish a correspondence, roughly speaking, between
Gorenstein perfect algebras of grade 1 that are birational onto their
image, on the one hand, and self-linked perfect ideals of grade 2 that
have one of the self-linking elements contained in the second symbolic
power, on the other hand.  Finally, we provide another characterization
of these ideals in terms of their symbolic Rees algebras, and we prove a
criterion for these algebras to be normal.
 \endabstract
 \endtopmatter
 \document

\newsect Introduction

A traditional way of studying an algebraic variety is to project it
from a general center onto a hypersurface, and conversely a variety
with given properties is often constructed by modifying a suitable
hypersurface.  It is known that the variety is Cohen--Macaulay if and
only if, locally, its algebra can be presented, as a module on the
target projective space, by a square matrix.  If so, then the
determinant cuts out the hypersurface, the submaximal minors generate
the `adjoint' ideal (the ideal that induces the conductor), and the
adjoint ideal is Cohen--Macaulay of codimension 2; moreover, if the
adjoint ideal is radical, then its second symbolic power contains the
determinant.  Furthermore, the variety is Gorenstein if and only if,
via row operations, the matrix can be made symmetric.  In this paper,
we'll prove in essence that the variety is Gorenstein if and only if
the adjoint ideal is self-linked with the determinant as one of the
self-linking elements; see Section~3.  In Section~2, we'll characterize
a self-linked Cohen--Macaulay ideal of codimension 2 as one with a
Hilbert--Burch matrix having a maximal symmetric subblock.  This result
was known before, but not in characteristic 2; moreover, our proof
yields, at the same time, the symmetrization result for Gorenstein
varieties.  In Section~3, we'll also characterize, in terms of the
symbolic Rees algebra, each such ideal that, in addition, has one of
the self-linking elements in its second symbolic power.  Finally, in
Section~4, we'll study what it means for the variety to be normal.

A milestone in this theory was reached by Catanese \cite{Cat} on the
basis of work by and with a number of others.  Catanese (see \cite{Cat,
Rmk.~5.1, \p.98}) set up a correspondence between Gorenstein algebras
with module generators and symmetric matrices satisfying a {\it rank
condition,} or {\it row condition,} (RC); here each algebra is
birational onto its image, the generators form a finite minimal set
starting with 1, and RC means that the ideal of submaximal minors of
the matrix is equal to the ideal of maximal minors of the submatrix
obtained by deleting the first row.  Viewing RC as the hypothesis of
the Rouch\'e--Capelli theorem (which asserts that an inhomogeneous
system of linear equations has a solution if the rank of the
coefficient matrix is equal to that of its augmentation), Catanese
\cite{Cat, Rmk.~4.6, \p.91} showed that RC was essentially equivalent
to the existence of a commutative associative algebra structure on the
cokernel of the matrix.  Catanese also considered algebras of higher
rank, especially of rank two, and represented them via symmetric
matrices; however, RC is no longer relevant.  (In fact, Catanese
considered weighted homogeneous algebras, but the theory is essentially
the same for our purposes.)  Reviewing Catanese's paper, Reid [MR
86c:14027] called the correspondence ``the most interesting idea of
this important paper.''

The correspondence was extended by Mond and Pellikaan \cite{MP} and by
de Jong and van Straten \cite{dJvS}.  Mond and Pellikaan were
especially interested in the locus $N_i$, for $i=1,2,3,\dots$, defined
by the vanishing of the $(i-1)$-st Fitting ideal of the algebra.  Their
purpose was to show that these loci are suitable choices for the loci
of multiple points of the projection.  For example, $N_1$ is the image of
the projection; $N_2$ is the locus where the projection is not an
isomorphism, with the structure defined by the adjoint ideal; and $N_3$
is, in the Gorenstein case, a particularly well-behaved locus of triple
points.  De Jong and van Straten were especially interested in relating
the deformation theory of the projection to the deformation theory of
the inclusion of $N_2$ into $N_1$.  In particular, they often worked
over a nonreduced base ring containing the complex numbers.

Mond and Pellikaan and de Jong and van Straten extended the
correspondence to relate Cohen--Macaulay algebras with generators and
square matrices with RC.  Furthermore, they added a third component,
the Cohen--Macaulay ideals $I$ of codimension 2 with a preferred
element $\D$ satisfying a {\it ring condition\/}; this is the condition
that the reciprocal of $I/(\D)$, viewed as a fractional ideal, is equal
to the endomorphism ring of $I/(\D)$.  The present authors, together
with Lipman, developed the connection with enumerative multiple-point
theory in \cite{KLU92} and \cite{KLU}.  In the present paper, we will
investigate the third component of the correspondence, especially in
the Gorenstein case in part using the symbolic Rees algebra, and we'll
study the significance of normality.  Along the way, we'll extend some
of the known results, and simplify and clarify their proofs.  In more
detail, here's what we'll do.

We will work, for the most part, over an arbitrary Noetherian local ring
$R$, which plays the role of the localized polynomial ring.  In this
setting, the correspondence relates $R$-algebras $B$ with module
generators, square matrices $\varphi$ with RC, and ideals $I$ with a
preferred element $\D$ satisfying the ring condition.  Moreover, we'll
replace the condition that $R$ be regular and $B$ (resp., $I$) be
Cohen--Macaulay with the more general condition that $B$ be perfect of
grade 1 (resp., that $I$ be perfect of grade 2).  Similarly, we'll use
the relative form of the Gorenstein condition, namely, that
$\Ext_R^1(B,R)$ be $B$-isomorphic to $B$.

In Section~2, we will prove the symmetrization results.  First,
Theorem~(2.3) asserts that, if a Gorenstein $R$-algebra $B$ is presented
by a matrix $\varphi$ with regular determinant, then there exist
invertible matrices $\varepsilon$ and $\mu$ such that
$\varepsilon\varphi$ and $\varphi\mu$ are symmetric.  Conversely, an
$R$-algebra $B$ is Gorenstein if it's presented by a matrix with regular
determinant $\D$ and if $R/(\D)\subset B$; see Proposition~(2.12).
Second, Theorem~(2.6) asserts that a necessary and sufficient condition
for a perfect $R$-ideal $I$ of grade two to be self-linked is that,
given any $n$ by $n-1$ matrix $\phi$ presenting $I$, there exists an
invertible matrix $\varepsilon$ such that the submatrix consisting of
the last $n-1$ rows of $\varepsilon\phi$ is symmetric.  Both
Theorem~(2.3) and necessity in Theorem~(2.6) are derived from
Lemma~(2.1), which shifts the burden of proving the symmetry of a matrix
to checking the condition that the syzygy matrix of the matrix modulo a
nonzerodivisor is symmetric, which is a natural condition in the context
of the two theorems.  Lemma~(2.1) enters via Proposition~(2.2), which
adds another common ingredient, a bilinear form.  This line of reasoning
was inspired by Catanese's original argument \cite{Cat, \pp.84--87}.
Theorem~(2.3) was proved independently and in essentially the same way
by Grassi \cite{Grassi}, and he went on to treat perfect Gorenstein
algebras of grade 2.  Another proof of Theorem~(2.3) was given by Mond
and Pellikaan \cite{MP, Prop.~2.5, \p.117}, but their proof uses Noether
normalization, and breaks down when $B$ is not finite over a regular
subring.

The sufficiency assertion of Theorem~(2.6) is derived in a few lines
from Part~(3) of Lemma~(2.5), which asserts this: given a symmetric
matrix whose ideal of submaximal minors that do not involve the first
column is of grade at least 2, the determinant and the lower right
minor together form a regular sequence with respect to which the ideal
is self-linked.  This assertion also follows from a theorem of Valla's
\cite{Valla84, Thm.~2.1, \p.97}; Valla's proof uses the mapping cone,
whereas ours uses matrix factorization.  The use of matrix
factorization in this context was introduced by de Jong and van Straten
\cite{dJvS, \p.532}, who used it to prove a version of Parts~(1) and
(2) of Lemma~(2.5); see also \cite{MP, pf.~of 3.14, \pp.126--7}.
Theorem~(2.6) was first stated and proved in full by Ferrand in an
unpublished fifteen-page manuscript, but his proof of necessity, unlike
ours, requires that 2 be invertible in $R$.  Valla \cite{Valla84,
Remark, \p.99} posed the problem of finding an easier proof of
necessity, and then gave one in the case of three generators.

Proposition~(2.9) gives, for a perfect $R$-algebra $B$ of grade 1, an
upper bound on the height of its Fitting ideal $F_i$ for $i\ge2$.  In
geometric terms, this is an upper bound on the codimension of the
multiple-point locus $N_{i+1}$ of a central projection.  Part~(1) of the
proposition asserts that the height of $F_i$ is at most ${i+1\choose2}$
if $B$ is Gorenstein.  Part~(2) asserts that the height of $F_2$ is at
most 4 if $B$ is Cohen--Macaulay and is of rank 2 over $R/(\D)$, where
$\D$ is a nonzerodivisor, and if $R$ is regular.  Part~(1) for $i=2$ was
proved (implicitly) by Mond and Pellikaan \cite{MP, Thm.~4.3, \p.131} when
$B$ is of rank 1 over its image (or birational onto it).  In particular,
in the Gorenstein case, the codimension of $N_3$ is at most 3, which is
the expected upper bound.  However, this bound can fail in the
Cohen--Macaulay case; see Example~(2.10).
 % for instance, let $I$ be the ideal of maximal
 % minors of an $n$ by $n-1$ matrix of indeterminates where $n\ge3$, let
 % $\D\in I^2$, and let $B$ be as in Proposition~(3.1).
 It is unknown whether it is possible to have $\cod N_4>4$ in the
Gorenstein case.

Both parts of Proposition~(2.9) will be derived from Lemma~(2.8), which
asserts that $F_i$ has the same radical as a certain subideal, the
subideal of $n-i$ by $n-i$ minors of the matrix obtained by deleting the
first row and column of a suitable matrix presenting $B$.  In Part~(1)
the matrix is symmetric, and in Part~(2) it's alternating.  In their
case, Mond and Pellikaan \cite{MP, Prop.~4.1, \p.128} proved that the
two ideals coincide.  We'll recover this result, via their approach, as
part of Lemma~(3.6).  In Part~(2) the matrix is given by
Proposition~(2.7)(1), which deals more generally with modules $M$ of
rank 2 over $R/(\D)$; the matrix is obtained directly from
Proposition~(2.2).  Proposition~(2.7)(2) recovers a result of Herzog and
K\"uhl \cite{HeK, Thm.~(3.1)(a), \p.82}, which asserts that the minimal
number of generators $\nu(M)$ is even.  They follow a different
approach: they associate a Gorenstein ideal of grade 3 to $M$ via a
Bourbaki sequence, and then use the Buchsbaum--Eisenbud structure
theorem, which says that these ideals are presented by an alternating
matrix.  Proposition~(2.7)(3) gives an upper bounds on the height of
$\Fitt_i(M)$ for $i<\nu(M)/2$; the special case $i=1$ is also a special
case of a result of Bruns \cite{Bruns, Cor.~2, \p.23}.

In Section~3, we'll establish the correspondence between the set of
perfect $R$-algebras $B$ of grade 1 whose first and second Fitting
ideals $F_1$ and $F_2$ satisfy $\grade F_1=2$ and $\grade F_2\ge3$ and
the set of the perfect $R$-ideals $I$ of grade 2 such that $I$ is a
complete intersection at each associated prime and whose second
symbolic power $\I2$ contains a preferred regular element $\D$,
determined up to a unit multiple.  The condition $\grade F_1=2$ means
simply that $B$ lies in the total quotient ring of $R/(\D)$ and that
$B\not=R/(\D)$; see, Proposition~(2.12).  The Gorenstein $B$ correspond
to the $I$ that are self-linked with $\D$ as one of the self-linking
elements.  Given $B$, we'll take $I$ to be $F_1$, and $\D$ to be a
generator of $F_0$.  Given $I$ and $\D$, we'll take $B$ to be the
reciprocal of $I/(\D)$ in the total ring of quotients of $R/(\D)$.
  Proposition~(3.1) starts with an $I$ and $\D$; it asserts that the
corresponding $B$ has the desired properties~--- in particular, that
$B$ is a ring, whose first Fitting ideal $F_1$ is equal to $I$, and
that $B$ is Gorenstein if and only if $I=(\D,\alpha):I$ for some
$\alpha$.  The proof shows that $B$ is presented by a square matrix
with RC.  In fact, its transpose presents $I/(\D)$ and the transpose of
its truncation is a Hilbert--Burch matrix of $I$.  Hence $F_1$ is equal
to the first Fitting ideal of $I/(\D)$, and RC means simply that
$F_1=I$.
  Conversely, Proposition~(3.6) starts with an algebra $B$ and asserts
that its first Fitting ideal has most of the desired properties; the
remaining ones follow from Proposition~(3.1).  The correspondence is
summarized formally in Theorem~(3.7).  All these results and their
proofs are found, by and large, in Mond and Pellikaan \cite{MP} and de
Jong and van Straten \cite{dJvS}.

This version of the correspondence features two conditions: (i)~that
$\D\in\I2$, and (ii)~that $I$ is self-linked with $\D$ as one of the
self-linking elements.  In the geometric case, under the assumption
that $R/I$ is reduced, the equivalence of (i) with the ring condition
was explicitly proved by de Jong and van Straten \cite{dJvS,
Thm.~(1.2), \p.529}.  In this situation, $I$ is generically a complete
intersection; so their case is covered by ours.  Moreover, in their
case, the map $d\:I/I^2\to\Omega^1_R\otimes R/I$ is generically
injective; hence, (i) means, in geometric terms, that the hypersurface
defined by the vanishing of $\D$ is singular along the zero locus of
$I$.  On a different tack, Eisenbud and Mazur, inspired by the recent
Taylor--Wiles proof of Fermat's last theorem, proved the following
higher-codimensional version of the equivalence of (i) and RC
\cite{Ei-Ma, Thm.~4}: if $I$ has no embedded primes and has grade at
least $r$, then (i) implies that the $(r-1)$-st Fitting ideal of
$I/(\D)$ is contained in $I$, and the converse holds if also $I$ is a
complete intersection of height $r$ locally at each associated prime.
As to (ii), it is equivalent to the condition that $B$ is Gorenstein,
given that $B$ is a ring.  Indeed, it is not hard to see that $\alpha$
is simply an adjoint; that is, $\alpha B=I/(\D)$.  Nevertheless, this
equivalence seems not to have been observed before.

In addition to having intrinsic interest, these conditions (i) and (ii)
suggest considering the symbolic Rees algebra
$\SR(I):=\bigoplus_{n\ge0}I^{(n)}t^n$.  Indeed, pursuing the work of a
number of others, Herzog and Ulrich \cite{H-U, Cor.~2.5, \p.146} proved
that, if $\SR I=R[It,\I2 t^2]$, then $I$ is self-linked in the case
that $R$ is regular of dimension 3 and that $I$ is of codimension 2 and
satisfies certain technical conditions.  We will, in essence,
generalize this result.  Indeed, Theorem~(3.4) asserts that, if $R$ is
Cohen--Macaulay with infinite residue class field and if $I$ is perfect
of grade $2$ and is a complete intersection at each associated prime,
then a necessary and sufficient condition for $I$ to be self-linked
with one of the self-linking elements in $\I2$ is that $\SR I=R[It,\D
t^2]$ for some $\D\in\I2$ and that $\I i$ be perfect for every $i>0$.
In fact, necessity holds without the Cohen--Macaulay hypothesis, see
Proposition~(3.3), and we suspect that sufficiency does too.  Moreover,
if $I$ is self-linked with respect to some $\D\in\I2$ and $\alpha\in
I$, then any general $\D$ and $\alpha$ will do; this is an unusual
phenomenon, and means that, in examples, a random choice of $\D$ and
$\alpha$ will likely do.  Theorem~(3.4) also gives another necessary
and sufficient condition, namely, that $\I2/I^2$ be cyclic and the
analytic spread $\ell(\I 2)$ be 2.

The quotient $\I2/I^2$ is rather interesting.  Indeed, Proposition~(3.8)
asserts that, if $B$ is Gorenstein, then $\Delta$ defines an
isomorphism, $R/F_2 \risom I^{(2)}/I^2$.  This result was proved under
the additional hypothesis that $F_1$ is radical by Mond and Pellikaan
\cite{MP, Thm.~4.4, \p.132}.  In addition, they pointed out its
``somewhat surprising'' geometric significance: the triple-point locus
of a projection depends only on the double-point locus, provided these
loci have at least the expected codimensions.  Another proof was given
by de Jong and van Straten \cite{dJvS, Thm.~(2.8), \p.538}.

In Section~4, we'll study the normality of $B$ when $R$ satisfies
Serre's conditions (R$_2$) and (S$_3$).  Proposition~(4.1) asserts
notably that, if $B$ is normal, then $\D\in\I2$ if and only if
$A:=R/(\D)$ has multiplicity 2 locally at each associated prime of $I$.
Consequently, it is easy to find an example in which $B$ is a normal
ring, yet $\D\notin\I2$; see Example~(4.2).  Finally, Theorem~(4.4)
asserts notably that, if $R$ contains an infinite and perfect field and
if $\D$ is a general element of $\I2$, then $B$ is normal if and only if
$A$ has multiplicity 2 locally at each associated prime of $I$;
furthermore, if so, then $B/A$ is unramified in codimension 1.  In other
words, if $\D\in\I2$ is general, then the hypersurface it defines is the
image of the projection of a normal variety with $I$ as its adjoint ideal
if and only if the hypersurface, generically along the zeros of $I$,
consists of two smooth sheets, which may be tangent.  Example~(4.5)
shows that it is necessary that $\D$ be taken general.

\newsect Symmetric Matrices

In this section, we prove that, if a Gorenstein algebra $B$ can be
presented by a square matrix $\varphi$ with regular determinant, then we
can take $\varphi$ to be symmetric; in fact, we can make $\varphi$
symmetric via row operations, or if we prefer via column operations.  We
use this result to derive a bound on the heights of the Fitting ideals
of $B$.  We also prove that a perfect ideal $I$ of grade 2 is
self-linked if and only if, via row operations, we can make any given
$n$ by $n-1$ matrix presenting $I$ possess a symmetric $n-1$ by $n-1$
subblock.

\lem1
 Let $(R,\m )$ be a  local ring, and let lower-case Greek
letters stand for $n$ by $n$ matrices with entries in $R$.  Let $\varphi$
be such a matrix, and let $r$ denote its rank modulo $\m $.  Let
$\D$ be a nonzerodivisor in $R$ contained in the annihilator
of  $\Cok\varphi$, and let `$\,\?{\phantom{s}}$' indicate
reduction modulo $\D$.  Further, assume that there is an exact sequence,
	$$\?R^n @>\?\psi>> \?R^n @>{\?\varphi}>> \?R^n\dno1.1$$
 where $\?\psi$ is either symmetric or alternating; in the latter
case, assume that $r$ is even.  Then there exist invertible matrices
$\varepsilon$ and $\mu$ such that the products
$\varepsilon\varphi$ and $\varphi\mu$ are symmetric or
alternating, respectively.
 \pf
 It suffices to show that there exist invertible matrices $\alpha$,
$\beta$ such that $\alpha\varphi\beta$ is symmetric (or alternating).
Indeed, let `$^*$' indicate transpose, and set
$\varepsilon:=\beta^{-1*}\alpha$ and $\mu:=\beta\alpha^{-1*}$.  Then
$\varepsilon\varphi$ is equal to
$\beta^{-1*}\alpha\varphi\beta\beta^{-1}$, which is obviously then
symmetric (or alternating).  Similarly, $\varphi\mu$ is symmetric
(or alternating).

Because $(R,\m )$ is  local, there exist invertible
matrices $\alpha$, $\beta$ such that $\alpha\varphi\beta$ has the form,
	$$\pmatrix \bold 1_{r\x r} &\bold 0 \\
		   \bold 0	   &\varphi'\endpmatrix,$$
 or, if $r=2s$, the form,
      $$\pmatrix \bold 0	 &\bold 1_{s\x s}&\bold 0 \\
		 -\bold 1_{s\x s}&\bold 0	 &\bold 0 \\
		 \bold 0	 &\bold 0	 &\varphi'\\
	\endpmatrix,$$
 where the entries of $\varphi'$ are in $\m $.  Replacing $\varphi$
by $\alpha\varphi\beta$ and $\psi$ by $\beta^{-1}\psi\beta^{-1*}$, we
may assume that $\varphi$ itself has the displayed form.  Now,
$\?\psi\?\varphi=\bold 0$, and $\?\psi$ is symmetric (or alternating).
Hence, we may assume that $\?\psi$ has the form,
	$$\pmatrix \?{\bold 0}_{r\x r} &\?{\bold 0} \\
		   \?{\bold 0}	     &\?\psi'\endpmatrix,$$
 where $\?\psi'$ is still symmetric (or alternating).  Therefore, the
following sequence is exact:
   $$\?R^{n-r} @>{\?\psi}'>> \?R^{n-r} @>{\?\varphi}'>> \?R^{n-r}.$$
 So, replacing $\varphi$ by $\varphi'$, we may assume that all the
entries of $\varphi$ lie in $\m $.

We may further assume that $\psi$ is symmetric (or alternating).  Using
the exact sequence \Cs1.1), we now prove (in a manner reminiscent of an
argument used by Shamash \cite{Shamash, pf. of Lem.~1, \p.454}) and by
  Eisenbud \cite{Eisenbud, 6.3,\p.55})
 that the inverse matrix $\psi^{-1}$ is defined over the total quotient
ring $K$ of $R$ and that there exists an invertible matrix $\varepsilon$
such that $\varepsilon\varphi= \D\psi^{-1}$.

First of all, $\D\cdot\bold 1_{n\x n}=\varphi\chi$ for some $\chi$;
indeed, since $\D$ annihilates the cokernel of $\varphi$, every column
of $\D\cdot\bold 1_{n\x n}$ is a relation on $\Cok\varphi$, and so lies
in the column space of $\varphi$.  Since $\D$ is a nonzerodivisor,
$\varphi^{-1}$ and $\chi^{-1}$ exist over $K$ and
	$$\chi=\D\varphi^{-1}.\dno1.2$$
   Consequently,
	$$\chi\varphi=\varphi\chi=\D\cdot\bold 1_{n\x n}.\dno1.3$$
 Since $\?\varphi\?\chi$ vanishes and since \Cs1.1) is exact, there
exists a $\gamma$ such that $\?\chi=\?\psi\?\gamma$.  Hence, for some
$\lambda$,
	$$\chi=\psi\gamma+\D\lambda=\psi\gamma+\chi\varphi\lambda,$$
 where the second equality follows from \Cs1.3).  Now, the entries of
$\varphi$ lie in $m$; hence, $1-\varphi\lambda$ is invertible.  Set
$\varepsilon:=\gamma(1-\varphi\lambda)^{-1}$.  Then
	$$\chi=\psi\varepsilon.\dno1.4$$
 Hence, since $\chi$ is invertible over $K$, so are $\psi$ and
$\varepsilon$.

On the other hand, since \Cs1.1) is a complex, there exists some $\eta$
such that
	$$\varphi\psi=\D\eta=\varphi\chi\eta.$$
 Multiplication by $\varphi^{-1}$ now yields $\psi=\chi\eta$.
Substituting into \Cs1.4), we conclude that $\varepsilon$ is, in fact,
invertible over $R$.  Finally, \Cs1.4) and \Cs1.2) yield that
$\varepsilon\varphi=\D\psi^{-1}$.  Now, $\psi^{-1}$ is symmetric (or
alternating).  Hence $\alpha:=\varepsilon$ and $\beta:=\bold 1_{n\x n}$
are as required, and the proof is complete.

\prp2
 Let $(R,\m )$ be a  local ring.  Let $\D$ be a nonzerodivisor
in $R$, and set $A:=R/(\D)$.  Let $M$ be an $A$-module that is
presented, as an $R$-module, by an $n$ by $n$ matrix $\varphi$.  Assume
that there exists a bilinear form $b\:M\x M\to A$ which induces a
surjection,
			 $$M\onto\Hom_A(M,A),$$
 and which is either symmetric or alternating.  In the latter case,
assume that
	 $$n\equiv\nu(M)\mod2,$$
 where $\nu(M)$ denotes the minimal number of generators of $M$.  Then
there are invertible $n$ by $n$ matrices $\varepsilon$ and $\mu$ such
that the products $\varepsilon\varphi$ and $\varphi\mu$ are
symmetric or alternating, respectively.
 \pf
 Let $u_1$, \dots, $u_n$ be generators of $M$ corresponding to the
presentation $\varphi$.  Let `$\,\?{\phantom{s}}$' indicate images in
$A$, and let `$^*$' indicate the dual of a given $A$-module or $A$-map.
Now, dualizing the exact sequence,
	$$A^n @>\?\varphi>>A^n@>>>M@>>>0,$$
 yields an exact sequence,
	$$0@>>>M^*@>>>A^{n*} @>\?\varphi^*>>A^{n*}.\dno2.1$$
 Set $\?\psi:=\bigl(b(u_i,u_j)\bigr)$.  This is a symmetric (or
alternating) $n$ by $n$ matrix with entries in $A$.  Since the $u_i$
correspond to the presentation of $M$ by $\varphi$, our assumption on
$b$ and the exactness of \Cs2.1) imply the exactness of the following
sequence:
	$$A^n@>\?\psi>>A^{n*} @>\?\varphi^*>>A^{n*}.$$ Obviously, modulo
$\m$, the rank of $\varphi$ is equal to $n-\nu(M)$.  Hence, Lemma~\Cs1)
yields the assertion.

\thm3
 Let $R$ be a Noetherian local ring, $B$ a Gorenstein $R$-algebra.  If
$B$ is presented, as an $R$-module, by an $n$ by $n$ matrix $\varphi$
with regular determinant, then there exist invertible $n$ by $n$
matrices $\varepsilon$ and $\mu$ such that $\varepsilon\varphi$ and
$\varphi\mu$ are symmetric.
 \pf
 Set $\D:=\det\varphi$ and $A:=R/(\D)$.  Since $\D$ is $R$-regular and
kills $B$, the sequence
	$0\to R @>\D>> R\to A\to0$
 is exact and yields a natural isomorphism of $B$-modules, $\hom_A(B,A)
= \ext1R(B,R)$.  This Ext is $B$-isomorphic to $B$  since $B$ is
Gorenstein over $R$.  Therefore, for some
$t\in\hom_A(B,A)$,
			  $$\hom_A(B,A)=Bt.\dno3.1$$

Let $m\:B\x B\to B$ be multiplication on the commutative ring $B$, and
consider the symmetric $A$-valued bilinear form $b:=t\circ m$.  The
equality \Cs3.1) implies that $b$ induces a surjection,
			 $$B\onto\Hom_A(B,A).$$
 Hence, the assertion follows from Proposition \Cs2).

\proclaim Lemma{ \(Laplace identity)} 4 Let $R$ be a commutative ring,
$\varphi$ an $n$ by $n$ matrix with entries in $R$.  Set
$\D:=\det\varphi$, denote by $m_i^j$ the minor of $\varphi$ obtained by
deleting Row $i$ and Column $j$, and denote by $m_{i,k}^{j,l}$ that
obtained by deleting Rows $i$ and $k$ and Columns $j$ and $l$, where
$m_{i,k}^{j,l}=0$ if $i=k$ or $j=l$, and where the determinant of the
empty matrix is $1$.  Then, for $1\le i\le k\le n$ and $1\le j\le
l\le n$,
	$$m_{i,k}^{j,l}\D=m_i^jm_k^l-m_i^lm_k^j.$$
 \pf We may assume that $1=i<k$ and $1=j<l$.  Let $\phi'$ be the
$n-1$ by $n-1$ matrix obtained by deleting the first row and column of
$\varphi$, and let $\theta$ be the $n$ by $n$ matrix obtained from
$\adj(\phi')$ by adding a row of 0s at the top and a column of 0s
on the left.  We can now check the asserted identity by multiplying the
familiar equation,
	$$\D\bold 1_{n\x n}=\varphi\adj(\varphi),$$
 on the left by $\theta$.

\lem5 Let $R$ be a Noetherian local ring, $\varphi$ an $n$ by $n$
matrix with entries in $R$.  Let $\phi$ be the $n$ by $n-1$ matrix
consisting of the last $n-1$ columns of $\varphi$, and $\phi'$ the
submatrix of $\phi$ consisting of the last $n-1$ rows.  Set
$\D:=\det\varphi$, set $\alpha:=\det\phi'$, and set $I:=I_{n-1}(\phi)$,
the ideal of $n-1$ by $n-1$ minors.  Assume $\grade I\ge2$.  Set
$A:=R/(\D)$ and let `$\,\?{\phantom{s}}\,$' indicate the image in $A$.
 \part 1 Let `$\null\,{}^*$' indicate the transpose.  If $\D$ is regular, then
 $$\gather \?I=\Cok(\?\varphi)=\Im(\adj\?\varphi)=\Ker(\?\varphi),\and\\
	\Hom_A(\?I,A)=\Ker(\?\varphi^*)
		=\Im(\adj\?\varphi^*)=\Cok(\?\varphi^*).\endgather$$
 \part 2 If $\D$ is regular, then
 $$\BI_{n-1}(\varphi)=\BI_{n-1}(\phi) \text{ if and only if }
	\Hom_A(\?I,A)=\Hom_A(\?I,\?I).$$
 \part 3 If $\varphi$ is symmetric, then $\D,\alpha$ form an
$R$-regular sequence, and $I$ is self-linked with respect to
$\D,\alpha$; that is, $(\D,\alpha):I=I$.
 \pf
 Let $\D_1$, \dots, $\D_n$ be the signed maximal minors of $\phi$, and
say
	$$\varphi=\pmatrix a_1\\
			    \vdots&\quad\up\phi\quad\\
			    a_n\endpmatrix.$$
 Then $\D=\sum a_i\D_i$.  Since $\grade I\ge2$, the Hilbert--Burch sequence,
	$$R^{n-1} @>{\phi}>> R^n @>>> I @>>> 0,$$
 is exact.  Hence, it induces this exact sequence:
	$$A^n @>{\?\varphi}>> A^n @>>> \?I @>>> 0.$$
 Dualizing the latter yields the following exact sequence:
	$$0\to \Hom_A(\?I,A)\to A^{n*}@>{\?\varphi}^*>>A^{n*}.$$
 On the other hand, it is not hard to show that, if $\D$ is regular,
then the matrix factorization $\varphi\cdot\adj\varphi=\D\bold 1_{n\x
n}$ yields the following two periodic exact sequences (see
\cite{Eisenbud, Prop.~5.1, \p.49}):
	\def\barvarphi{\?\varphi}
     $$\gather\cdots@>>> A^{n} @>\barvarphi>> A^{n} @>\adj\barvarphi>>
		A^{n} @>\barvarphi>> A^{n} @>>>\cdots;\\
	\cdots@>>> A^{n*} @>\barvarphi^*>> A^{n*} @>\adj\barvarphi^*>>
		A^{n*} @>\barvarphi^*>> A^{n*} @>>>\cdots.\endgather$$
 Therefore, (1) holds.

Assume that $\D$ is regular.  Let $u_i\:\?I\to A$, where $1\le i\le n$,
be the generators of $\Hom_A(\?I,A)$ that arise, in the above
presentation of this module, from the dual basis of $A^{n*}$.  Then
	$$\adj\barvarphi^* = \bigl(u_i(\?\D_j)\bigr).$$
 Thus $\BI_{n-1}(\?\varphi)\?{}=\?I$ if and only if
$u_i(\?I)\subset\?I$ for every $i$.  It follows that (2) holds.

Assertion (3) holds by virtue of some of Valla's work; see
\cite{Valla84, Thm.~2.1, \p.97}.  However, the following proof is more
in the spirit of the present paper.  Use the notation of the preceding
lemma.  Then $m_1^1=\alpha$ and $m_i^1=(-1)^{i+1}\D_i$.  So the
preceding lemma (the Laplace identity) yields
	$$m_{1,k}^{j,1}\D=m_1^jm_k^1-\alpha m_k^j.$$
 Assume $\varphi$ is symmetric.  Then $m_1^j=m_j^1$.  Hence
$(\D,\alpha)\supset I^2$.  Consequently, $\D,\alpha$ form an
$R$-regular sequence.

Let $K$ denote the total ring of quotients of $A$.  Then $\?IK=K$
because $\?I$ contains a nonzerodivisor, namely, $\?\alpha$.  View
$\Hom_A(\?I,A)$ as a submodule of $\Hom_K(\?IK,K)$, and view the latter
as $K$ by identifying a map $u$ with multiplication by $u(1)$.  Then
$\Hom_A(\?I,A)$ becomes identified with $A:_K\?I$, and $\Hom_A(\?I,\?I)$
with $\?I:_K\?I$.  Consider the image of $\Hom_A(\?I,A)$ in $A^{n*}$,
and project it on the first coordinate; clearly, the projection is
$(A:_K\?I)\?\alpha$.  The latter is equal to $(\?\alpha):_K\?I$, so to
$(\?\alpha):_A\?I$ because $\?\alpha\in\?I$ is a nonzerodivisor.  On the
other hand, the image of $\Hom_A(\?I,A)$ is equal to
$\Im(\adj\?\varphi^*)$ by (1).  Hence its projection is equal to the
ideal generated by the first row of $\adj\?\varphi^*$.  Since $\varphi$
is symmetric, the latter is $(\?\D_1,\dots,\?\D_n)$, or $\?I$.  Thus
$(\?\alpha):\?I=\?I$, and so $(\D,\alpha):I=I$.  The proof is now
complete.

\thm6 Let $R$ be a Noetherian local ring, let $I$ be a perfect
$R$-ideal of grade two, and let $\phi$ be an $n$ by $n-1$ matrix with
entries in $R$ presenting $I$.  Then the following two conditions are
equivalent:
 \smallskip
 \item i  The ideal $I$ is self-linked.
 \item ii There exists an invertible $n$ by $n$ matrix $\varepsilon$
with entries in $R$ such that the $n-1$ by $n-1$ matrix consisting of
the last $n-1$ rows of $\varepsilon\phi$ is symmetric.
 \pf
 We first show that (i) implies (ii).  Assume that $(\D,\alpha):I=I$.
Then $\D,\alpha$ form an $R$-regular sequence contained in $I$, because
	$$I^2\subset(\D,\alpha)\subset(\D,\alpha):I=I.$$
  Let $\D_1$, \dots, $\D_n$ be the signed maximal minors of $\phi$.
They generate $I$ by the Hilbert--Burch Theorem, which applies because
$I$ is perfect of grade two.  So, since $\D$ is in $I$, there are
elements $a_1$, \dots, $a_n$ of $R$ such that $\D=\sum a_i\D_i$.  Set
	$$\varphi:=\pmatrix a_1\\
			    \vdots&\quad\up\phi\quad\\
			    a_n\endpmatrix.$$
 Then $\det\varphi=\D$.  Set $A:=R/(\D)$ and let
`$\,\?{\phantom{s}}\,$' indicate the image in $A$.  Notice that the
following sequence is exact:
	$$R^n @>{\varphi}>> R^n @>>> \?I @>>> 0.$$

Let $K$ denote the total ring of quotients of $A$, and $t:K\to K$
 multiplication by $1/\?\alpha$.  As in the proof of Lemma (2.5), we
can see that $\Hom_A(\?I,A)$ is equal to $A:_K\?I$.  The latter is
clearly equal to $\big(\?\alpha:_A\?I\bigr)(1/\?\alpha)$.  Hence, the
self-linkage assumption yields
	$$\Hom_A(\?I,A)=\?I t.\dno6.1$$

On the other hand, since $\?I^2\subset(\?\alpha)$, multiplication on
$A$, followed by $t$, defines a symmetric bilinear form $b\:\?I\x\?I\to
A$.  By virtue of (2.6.1), $b$ induces an isomorphism,
	$$\?I\risom\Hom_A(\?I,A).$$
 Proposition (2.2) now implies (ii).

Conversely, assume (ii).  Since $\varepsilon$ is invertible,
$\varepsilon\phi$ also presents $I$.  Hence, we may assume that
	$$\phi=\pmatrix a_2&\cdots &a_n\\
				&\hfill\phi'\hfill\endpmatrix$$
 for suitable elements $a_2$, \dots, $a_n$ and a suitable symmetric
$n-1$ by $n-1$ matrix $\phi'$.  Set
	$$\varphi:=\pmatrix 0&a_2&\cdots&a_n\\
			    a_2\\
			    \vdots&&\hfill\up{\phi'}\hfill\\
			    a_n\endpmatrix$$
 and apply Lemma (2.5).  Thus (ii) implies (i), and the proof is
complete.

\prp7 Let $R$ be a Noetherian local ring, $\D$ an $R$-regular element,
and set $A:=R/(\D)$.  Let $M$ be an $A$-module.  Assume that $M_p$ is a
free $A_p$-module of rank two for every prime ideal $p$ of $A$ with
$\depth A_p\le1$, that $M$ is orientable (that is,
${(\wedge^2M)}^{**}\cong A$ where $^*$ indicates the dual module), and
that $M$ is presented, as an $R$-module, by an $n$ by $n$ matrix
$\varphi$ with $n\equiv\nu(M)\mod2$, where $\nu(M)$ denotes the minimal
number of generators of $M$.  (For instance, let $M$ be an orientable
maximal Cohen--Macaulay module of rank two over a hypersurface ring.)
 \part 1 Then there are invertible $n$ by $n$ matrices $\varepsilon$
and $\mu$ with entries in $R$ such that $\varepsilon\varphi$ and
$\varphi\mu$ are alternating.
 \part 2 (Herzog and K\"uhl \cite{HeK, Thm.~(3.1)(a), \p.82}) Then $\nu(M)$ is
even.
 \part 3 For $i<\nu(M)/2$,
	$$\height\Fitt_{2i}^A(M)=\height\Fitt_{2i+1}^A(M)\le2i^2+3i.$$
 \pf To prove (1), consider the alternating bilinear form,
	$$b\:M\x M\to\wedge^2M\to(\wedge^2M)^{**}\risom A,$$
 where the first two maps are the natural maps; it induces a surjection
$M\to M^*$ since $M$ is free locally in depth 1.  Hence Proposition
(2.2) yields (1).

For (2), we may assume that $n=\nu(M)$ and that $\varphi$ is
alternating.  Notice that
	$$\D\in\Ann_R(M)\subset\sqrt{\Fitt_0^R(M)}.$$
 In particular, $\det\varphi\neq0$, and hence $n$ is even.

Consider (3).  (The special case $i=1$ is also a special case of
\cite{Bruns, Cor.~2, \p.23}, because $M$ is a second syzygy module of
itself over $A$.)  By (2), $n$ is even.  Hence, by \cite{BE77b,
Cor.2.6, \p.462},
	$$\sqrt{\Fitt_{2i}^R(M)}=\sqrt{\Fitt_{2i+1}^R(M)},$$
 and these two radicals are also equal to the radical of the ideal of
$n-2i$ Pfaffians of $\varphi$.  The latter ideal has height at most
$2i+2\choose2$ by \cite{JP, Thm.~2.1, \p.191}, where this statement is
shown to follow from \cite{EN, \pp.202--3}.  Since these radicals
contain $\D$, we conclude that
	$$\height\Fitt_{2i}^A(M)=\height\Fitt_{2i+1}^A(M)
	\le{2i+2\choose2}-1=2i^2+3i.$$

\lem8 Let $R$ be a local ring, $B$ an $R$-algebra.  Assume that, with
respect to a set of generators $1$, $u_2$, \dots, $u_n$, the $R$-module
$B$ can be presented by a symmetric or alternating $n$ by $n$ matrix
$\varphi$ with regular determinant.  Set $F_i:=\Fitt_i^R(B)$.  Let
$\phi'$ be the matrix obtained by deleting the first row and column of
$\varphi$, and set $F_i':=\BI_{n-i}(\phi')$.  Then, for $i\ge2$,
	$$\sqrt{F_i}=\sqrt{F_i'}.$$
 \pf Let $p$ be a prime ideal of $R$ containing $F_i'$.  We need to
show that $p$ contains $F_i$.  Localizing at $p$, we may assume that
$p$ is equal to the maximal ideal $m$ of $R$.  If some entry of $\phi'$
lies outside of $m$, then we may perform row and column operations to
reduce the size of $\varphi$, until $\phi'$  has all its entries in $m$.
Furthermore,  $F_i'\neq R$; so $n-i\ge1$.  Hence $n-1\ge i\ge2$.
Therefore, $F_0\subset m$, or equivalently, $B\neq0$.

Say $\varphi=(a_{ij})$.  We are going to show that the $a_{ij}$ lie in
$m$.  Suppose that $a_{1j}\notin m$ for some $j$ with $2\le j\le n$.
Then $1\in mB$ because the columns of $\varphi$ yield relations among
$1$, $u_2$, \dots, $u_n$ and because $\phi'$ has all its entries in
$m$.  Hence $B=mB$.  However, $B\neq0$.  Therefore, $a_{1j}\in m$ for
$2\le j\le n$, and so $a_{j1}\in m$ too because $\varphi$ is symmetric
or alternating.  Finally, repeating the argument with the first column
of $\varphi$, we see that $a_{11}\in m$ as well.  Thus $a_{ij}\in m$
for all $i,j$.  So $F_i\subset m$ since $n-i\ge1$.

\prop9 \(1) Let $R$ be a Noetherian local ring, $B$ a finite
$R$-algebra.  Assume that $B$ is a perfect $R$-module of grade $1$ and
that $B$ is Gorenstein over $R$.  Set $F_i:=\Fitt_i^R(B)$.  Then, for
$i\ge2$, either $\height F_i\le{i+1\choose2}$ or $F_i=R$.
 \part 2  Let $R$ be a regular local ring,  $\D$ an $R$-regular element,
and set $A:=R/(\D)$.  Let $B$  be an $A$-algebra.  Assume that $B$  is a
finitely generated $A$-module of rank $2$ and that $B$   is a
Cohen--Macaulay ring.  Set $F_2:=\Fitt_2^R(B)$.  Then either $\height
F_2\le4$ or $F_2=R$.
 \pf In (1), the assertion follows from Theorem~(2.3), from Lemma~(2.8),
and from \cite{Jzfk, Thm.~2.1, \p.597} applied to the $n-1$ by $n-1$ matrix
$\phi'$ of Lemma~(2.8).

Consider (2).  Suppose that $\height F_2 >4$; then
$\height\Fitt_2^A(B)>3$.  Hence, as an $A$-module, $B$ is free of rank
$2$ locally in codimension 3.  Therefore, the reflexive $A$-ideal
$(\wedge^2B)^{**}$ is principal locally in codimension 3.  Hence,
$(\wedge^2B)^{**}$ is globally principal by \cite{Groth, Thm.~3.13(ii),
\p.132}; in other words, $B$ is an orientable $A$-module.  Therefore,
Proposition~(2.7) applies, and the assertion follows from Lemma~(2.8) and
\cite{JP, Thm.~2.1, \p.191}.

\eg10 The bound in Proposition~(2.9)(1) can fail if $B$ isn't assumed
to be Gorenstein over $R$, even for $i=2$.  For example, let $n\ge3$,
let $R$ be the power series ring in $n(n-1)$ variables over a field,
and let $I$ be the ideal of maximal minors of an $n$ by $n-1$ matrix
$\phi$ in those variables.  Let $\D\in I^2$, and let $B$ be the
reciprocal of $I/(\D)$, viewed as a fractional ideal in the total ring
of quotients of $R/(\D)$.  Then $B$ is an $R$-algebra, which is a
perfect $R$-module of grade $1$, by Parts~(1) and (2) of
Proposition~(3.1) below.  Moreover, the proof of Part~(1) shows that
$F_2$ may be described as the ideal of $n-2$ by $n-2$ minors of the $n$
by $n$ matrix obtained by adding to $\phi$ the column of combining
coefficients in an expansion of $\D$ as a linear combination of the
maximal minors of $\phi$.  So $F_2$ contains the ideal of $n-2$ by
$n-2$ minors of $\phi$.  Therefore, $\height F_2\ge6$.  This example
also shows that Proposition~(2.9)(2) can fail if $B$ is of rank $1$.

\lem11
 Let $R$ be a Noetherian local ring.  Let $B$ be an $R$-algebra, and
assume that, as an $R$-module, $B$ is finitely generated and perfect of
grade $1$.  Set $A:=R/\Fitt_0^R(B)$, and let $K$ and $L$ be the total
rings of quotients of $A$ and $B$.  Then the following six conditions
are equivalent: \smallskip
 \settabs 2\columns
 \+\item i  $A\subset B$.	&\item ii $\grade(\Fitt_1^R(B))\ge2$.\cr
 \+\item iii $K\onto K\ox_AB$. &\item iv $K=K\ox_AB$.\cr
 \+\item v  $B\subset K$. &\item vi $K=L$.\cr
 If $R$ satisfies Serre's condition \Ser2, then the above six
conditions are equivalent to the following one:
 \item vii $f\:\Spec(B)\to\Spec(A)$ is birational onto its image.
 \pf
 Let $\varphi$ be an $n$ by $n$ matrix with entries in $R$ presenting
$B$, and set $\D:=\det\varphi$.  Then $\D$ is $R$-regular, and
$(\D)=\Fitt_0^R(B)$.

Conditions (i) and (ii) are equivalent.  Indeed, (i) says that
$(\D)=\Ann_R(B)$.  However,
	$$\Ann(B)=(\D):\Fitt_1^R(B)$$
   by \cite{BE77a, Theorem, \p.232}, which applies because $\D$ is an
$R$-regular element contained in $\Fitt_1^R(B)$.  Hence, (i) holds if
and only if $\Fitt_1^R(B)$ is not contained in any associated prime of
$(\D)$, and the latter condition is equivalent to (ii).

Conditions (ii) and (iii) are equivalent.  Indeed, (ii) says that, for
each  $p\in\Ass(A)$,
	$$\Fitt_1^A(B)=\Fitt_1^R(B)\cdot A\not\subset p.$$
 Equivalently, the $A_p$-module $A_p\ox_AB$  is cyclic.  Or,
equivalently again, the natural map from $A_p$ to $A_p\ox_AB$ is
surjective.  The latter condition is simply (iii).

Trivially, (iv) implies (iii), and the converse holds because (iii)
implies (i) and  $K$ is $A$-flat.

Condition (iv) implies (v) because, in any event, $B\subset K\ox_AB$.
Indeed, take $p\in\Ass_A(B)$ and let $q$ be the preimage of $p$ in $R$.
Since $B$ is a perfect $R$-module of grade one, $\depth R_q=1$ by the
Auslander--Buchsbaum formula.  On the other hand, $p$ contains the
$R$-regular element $\D$.  Hence,
	$$q\in\Ass_R\bigl(R/(\D)\bigr)=\Ass_R(A).$$
 Therefore, $p\in\Ass_A(A)$.  Consequently, every element of $A$
regular on $A$ is also regular on $B$, and so $B\subset K\ox_AB$.

Condition (v) implies (i) because $A\subset K$.  Condition (v) implies
(vi).  Indeed, if $B\subset K$, then every nonzerodivisor on $A$ is a
nonzerodivisor on $B$, and so $B\subset K\subset L$; whence, $K=L$
since $K$ is its own total ring of quotients.

The map $f$ is locally of flat dimension 1 by the equivalence of (iii)
and (vi) of \cite{KLU, (2.3)}.  Hence $f$ is locally of codimension 1
by (2.5)(1) of \cite{KLU, (2.5)}.  By \cite{KLU, (3.2)(1)}, therefore,
$f$ is birational onto its image if and only if $\Fitt_1^R(B)$ is not
contained in any minimal prime of $(\D)$.  If $R$ satisfies \Ser2, then
$(\D)$ has only minimal associated primes.  However, $\Fitt_1^R(B)$ is
not contained in any associated prime of $(\D)$ if and only if
Condition (ii) holds.

\prp12
 Let $R$ be a Noetherian local ring, and $B$ an $R$-algebra.  Assume
that, as an $R$-module, $B$ is finitely generated and perfect of grade
$1$.  Assume that one of the equivalent conditions (i)--(vi) of
Lemma~(2.11) holds.  Let $\varphi$ be an $n$ by $n$ matrix with entries
in $R$ presenting $B$.  Then there are invertible $n$ by $n$ matrices
$\varepsilon$ and $\mu$ with entries in $R$ such that
$\varepsilon\varphi\mu$ is symmetric if and only if $B$ is Gorenstein
over $R$.
 \pf
 If there is a symmetric matrix $\varphi$ such that the sequence,
		 $$0 \To R^n @>\varphi>>R^n\To B\To0,$$
 is exact, then there is an isomorphism of $A$-modules,
	$$\Ext_R^1(B,R)\cong B,$$
  and it is even $B$-linear since $A\subset B\subset K$, where $A$ and
$K$ are as in Lemma~\Cs11).  The converse holds by Theorem~(2.3).

\newsect Symbolic powers

In this section, we establish a bijective correspondence between the
set of perfect algebras $B$ of grade 1 and the set of perfect ideals
$I$ of grade 2 with a preferred member $\Delta$ of its second symbolic
power, determined up to unit multiple; the first (resp., the second)
Fitting ideal of $B$ must have grade at least 2 (resp., at least 3),
and $I$ must be a complete intersection locally at each associated
prime.  The Gorenstein $B$ correspond  to the self-linked $I$
with $\Delta$ as one of the self-linking elements.  Proposition~(3.1)
starts with an $I$ and $\D$, takes $B$ to be the reciprocal of
$I/(\D)$ in the total ring of quotients $K$ of $A:=R/(\D)$, and
shows that $B$ has the requisite properties.  Conversely,
Proposition~(3.6) starts with a $B$; takes $I$ to be its first Fitting
ideal, and shows that $I$ has most of the requisite properties; the
remaining ones follow from Proposition~(3.1).  We summarize the
correspondence formally in Theorem~(3.7).  In addition, we investigate
the symbolic Rees algebra $\SR I$; notably, Theorem~(3.4) gives a
criterion in terms of $\SR I$ for $I$ to be self-linked with $\Delta$
as one of the self-linking elements.

To recall the precise definitions, let $I$ be an ideal in a Noetherian
ring $R$, set $S:=R/I$, let $W$ be the complement in $R$ of the union
of all associated primes of $I$, and let $n\ge0$ be an integer.  Then
the n{\it th symbolic power\/} $I^{(n)}$ of $I$ is the preimage of
$I^nR_W$ in $R$.  For $n>0$, notice that $I^n\subset I^{(n)}\subset I$
and that $$(I^{(n)}\cap I^{n-1})/I^n$$ is the $S$-torsion of
$I^{n-1}/I^n$.  The {\it symbolic Rees algebra\/} of $I$ is the graded
subalgebra,
	 $$\SR I:=\bigoplus_{n\ge0}I^{(n)}t^n,$$
 of the polynomial ring $R[t]$ in one variable $t$.  This algebra need
not be Noetherian, even if $I$ happens to be a prime ideal in a regular
ring; see \cite{Roberts85}, \cite{Roberts90}, and \cite{G-N-W}.

\prp1 Let $R$ be a Noetherian local ring, $I$ a perfect ideal of grade
$2$, and $\D\in I$ a nonzerodivisor.  Set $A:=R/(\D)$, and let
`$\,\?{\phantom{s}}$' indicate the image in $A$.  Let $K$ be the total
ring of quotients of $A$, and set $B:=A:_K\?I$ and $F_i:=\Fitt_i^R(B)$.
 \part 1 Then $\?I$ and $B$ are perfect $R$-modules of grade $1$.  In
addition,
$$\Fitt_1^R(I)=I,\ F_0=(\D),\ F_i=\Fitt_i^R(\?I), \and \?I=A:_KB.$$
 \part 2 Then the following four conditions are equivalent: (a)~$B$ is
a ring; (b)~$B=\?I:_K\?I$; (c)~$B\?I=\?I$; and (d)~$F_1=I$.  Moreover,
those four conditions are implied by this fifth condition: (e)
$\D\in\I2$.
 \part $2'$ If $\grade F_2\ge 3$, then all five conditions (a) to (e) of\/
(2) are equivalent.
 \part 3
 If $I$ is a complete intersection at each associated prime, then
$\grade F_2\ge 3$, and the converse holds if one of the conditions (a)
to (e) of\/ (2) holds.
 \part 4 Assume that $B$ is a ring.  Let $\alpha\in I$.  Then the
following conditions are equivalent:
 (i)~$B$ is Gorenstein over $R$, and $\?\alpha\notin \?In$  for any
maximal ideal $n$ of $B$;
 (ii)~$\?\alpha$ is a nonzerodivisor and $B={\?\alpha\,}^{-1}\?I$;
 (iii)~$\?I=\?\alpha B$; (iv)~$\?I^2=\?\alpha\?I$;
 (v)~$\?I=(\?\alpha):_A\?I$; and (vi)~$I=(\D,\alpha):I$.
 \pf
 Let $\phi$ be an $n$ by $n-1$ matrix with entries in $R$ presenting
$I$.  Let $\D_1$, \dots, $\D_n$ be the signed maximal minors of $\phi$.
They generate $I$ by the Hilbert--Burch Theorem; in other word,
$\Fitt_1^R(I)=I$.  So, since $\D$ is in $I$, there are elements $a_1$,
\dots, $a_n$ of $R$ such that $\D=\sum a_i\D_i$.  Set
$$\varphi:=\pmatrix a_1\\
\vdots&\quad\up\phi\quad\\
a_n\endpmatrix.$$
 Then $\det\varphi=\D$.  Furthermore, $\varphi$ presents $\?I$.
Hence, Lemma~(2.5)(1) yields
$$\?I=\Ker(\?\varphi)\and\Hom_A(\?I,A)=\Cok(\?\varphi^*).$$
 Dualize the second equation and combine the result with the first; thus
$$\Hom_A(\Hom_A(\?I,A),A)=\?I.$$
 On the other hand, $\Cok(\?\varphi^*)=\Cok(\varphi^*)$ because $\D$
annihilates $\Cok(\varphi^*)$ since $\D=\det\varphi^*$.  Hence the
second equation above implies that $\Hom_A(\?I,A)$ is presented by
$\varphi^*$ and is perfect of grade 1.  Finally, $\?I$ contains a
nonzerodivisor because $\grade I\ge2$; hence, there is a natural
identification of $\Hom_A(\?I,A)$ with $B$, and of $\Hom_A(B,A)$ with
$A:_KB$; see the end of the proof of (2.5).  Therefore, (1) holds.

Consider (2).  Obviously, (c) implies (b), and (b) implies (a); also,
(a) implies (c) because $\?I=A:_KB$ by (1).  Now, it follows from the
preceding paragraph that
	$$F_1=\BI_{n-1}(\varphi),\ I=\BI_{n-1}(\phi),\
	\Hom_A(\?I,A)=B\and \Hom_A(\?I,\?I)=\?I:_K\?I.$$
 Hence, (d) and (b) are equivalent by Lemma~(2.5)(2).  Thus the four
conditions (a) to (d) are equivalent.

Finally, assume (e), that $\D\in\I2$.  We'll prove (d), that $F_1=I$.  Now,
    $$F_1=\BI_{n-1}(\varphi^*)\supset\BI_{n-1}(\phi)=I.$$
 Hence, it suffices to prove (d) after localizing at an arbitrary
associated prime of $I$.  But then $\I2=I^2$, so $\D\in I^2$.  Hence we
may assume that $a_i\in I$ because the $\D_i$ generate $I$.  Therefore,
$$F_1=\BI_{n-1}(\varphi^*)\subset\BI_{n-1}(\phi)+I=I.$$
 Hence $F_1=I$.  Thus (2) holds.

Consider ($2'$).  The Laplace identity (2.4) implies that
$F_2(\D)\subset F_1^2$.  Hence, if $F_1=I$, then $\D\in\I2$ because
$F_2$ lies in no associated prime of $I$ as $I$ is perfect of grade 2
and $\grade F_2\ge3$.

Consider (3).  First assume that $I$ is a complete intersection at each
associated prime.   Since $I$ is perfect of grade 2, every prime
$p$ such that $p\supset I$ and  $\depth R_p=2$ is associated to $I$.
Hence $I_p$ can be generated by two elements.  Therefore,
$\grade\Fitt_2^R(I)\ge3$.  Hence, $\grade F_2\ge 3$ because
    $$F_2=\BI_{n-2}(\varphi^*)\supset\BI_{n-2}(\phi)=\Fitt_2^R(I).$$

Conversely, assume that $\grade F_2\ge 3$ and that one of the conditions
(a) to (e) of (2) holds.  Then $\D\in\I2$ by ($2'$).  Now, let
$p\in\Ass(I)$, and replace $R$ by $R_p$.  Then $\I2=I^2$.  In addition,
$F_2=R$ because $\grade F_2\ge 3$, but $\depth R_p=2$.  Hence (1)
implies that $\?I$ is generated by two elements.  Lift them to $I$.
Those lifts generate $I$ because $\D\in\I2=I^2$.  Hence the lifts form a
regular sequence because $\grade I=2$.

Consider (4).  The canonical isomorphisms,
	$$\Ext_R^1(B,R)=\Hom_A(B,A)=A:_KB=\?I,$$
 are obviously $B$-linear.  Hence, $B$ is Gorenstein over $R$ if and
only if $\?I$ is $B$-isomorphic to $B$.  Therefore, (i) implies (iii),
and the converse holds because $\?I$ contains a nonzerodivisor.
Obviously, (ii) and (iii) are equivalent.  Obviously, (iii) implies
(iv).  Suppose (iv) holds.  Then $\?\alpha$ is a nonzerodivisor,
because $\?\alpha$ divides the square of any element of $\?I$, and
$\?I$ contains a nonzerodivisor.  So $({\?\alpha\,}^{-1}\?I)\?I=\?I$.
So ${\?\alpha\,}^{-1}\?I\subset B$ because $B=\?I:_K\?I$ by (2).
However, $\?\alpha B\subset \?I$ because $\alpha\in I$; so $B\subset
{\?\alpha\,}^{-1}\?I$.  So (ii) holds.  Thus (i) to (iv) are equivalent.

If (v) holds, then $\?\alpha$ is a nonzerodivisor, because, again,
$\?\alpha$ divides the square of any element of $\?I$.  Now, whenever
$\?\alpha$ is a nonzerodivisor, then, clearly,
   $$(\?\alpha):_A\?I=(\?\alpha):_K\?I=\?\alpha(A:_K\?I)=\?\alpha B.$$
 Hence, (iii) and (v) are equivalent.  Obviously, (v) and (vi) are
equivalent.

\lem2
 Let $R$ be a Noetherian local ring with infinite residue field, and let
$I$ be an ideal of grade $2$ that is a complete intersection at each
associated prime.  Assume that the analytic spread $\ell(\I2)$ of $\I2$
is equal to $2$, and that $\I2/I^2$ is a cyclic module.  Let $\D\in\I2$
be a general element, and $\alpha\in I$ an element whose image in
$I/(\D)$ is general.  Then $I=(\D,\alpha):I$.
 \pf
 Let `$\,\?{\phantom{s}}$' indicate the image in $R/(\D)$.  Now, since
$\D\in\I2$ is general, it is part of a set of $\ell(\I2)$ elements
generating a reduction $J$ of $\I2$.  Then $\?J$ is a reduction of
$\?{\I2}$, and so $\?J$ requires at least $\ell\big(\?{\I2}\big)$
generators.  Hence $\ell(\?{\I2}) \le \ell(\I2)-1$, and therefore
$\ell(\?{\I2})=1$.  On the other hand, $\I2=I^2+(\D)$; hence,
$\?{\I2}=\?{I^2}$.  Therefore, $\ell({\?I}^2)=1$.  Consequently,
$\ell(\?I)=1$.  Since $\?\alpha\in\?I$ is general, therefore $\?\alpha$
generates a reduction of $\?I$.  In particular, $\?I$ and $\?\alpha$
have the same radical.  So $I$ and $(\D,\alpha)$ have the same radical
too.  Therefore, every associated prime $p$ of $(\D,\alpha)$ is one of
$I$; indeed,  $\depth R_p=2$ and  $p\supset I$, whence
$\depth(R/I)_p=0$ since $I$ is perfect of grade 2.

It suffices to check the equation $I=(\D,\alpha):I$ locally at each (of
the finitely many) associated primes $p$ of $I$, since any associated
prime of $(\D,\alpha):I$ is one of $(\D,\alpha)$, so one of $I$.
Localizing at $p$, we may assume that $I$ is a complete intersection.
Let $\m$ denote the maximal ideal of $R$.  Then $\?\alpha\notin\?\m\?I$
because $\?\alpha\in \?I$ is general.  So $\alpha\notin\m I$.  Hence $I=
(\alpha, \beta)$ for some $\beta$, and $\alpha,\beta$ form a regular
sequence.  Hence the associated graded ring $\gr_I(R)$ is a polynomial
ring over $R/I$; see \cite{Mats, \p.125} for instance.  In particular,
$\gr_I(R)$ is torsion free over $R/I$; hence, $\I2=I^2$.  Therefore,
$\D\in I^2$.  So $\D=r\alpha^2+s\alpha\beta+t\beta^2$ for some $r,s,t$.
Now, $\sqrt{I}= \sqrt{(\D,\alpha)}$, so $I^i=(\D,\alpha)\cap I^i$ for
some $i$.  On the other hand, by the generality of $\D$ and $\alpha$,
their leading forms in $\gr_I(R)$ have degrees 2 and 1, and they form a
regular sequence; it follows that $(\D,\alpha)\cap I^i= \D I^{i-2}+
\alpha I^{i-1}$ (see \cite{VV,Thm. 2.3, \p.97} for instance).  Hence
$I^i= \D I^{i-2}+ \alpha I^{i-1}$.  Suppose $t\in\m$.  Then $I^i= \alpha
I^{i-1}+ tI^i$.  So Nakayama's lemma implies that $I^i=\alpha
I^{i-1}$. Hence $I^i\subset(\alpha)$.  However, this inclusion is
impossible because $I^i$ has grade 2.  Thus $t\notin\m$.  Hence
$(\D,\alpha)=(\alpha,\beta^2)$.  Hence,
  $$(\D,\alpha):I=(\alpha,\beta^2):(\alpha,\beta)=(\alpha,\beta)=I.$$

\prp3 Let $R$ be a Noetherian local ring, $I$ a perfect ideal of grade
$2$.  Assume that $I$ is a complete intersection at each associated
prime, and that $I=(\D,\alpha):I$ for some $\D\in\I2$ and some
$\alpha\in I$.
 \part 1 If also $R$ has an infinite residue field, then $I=
(\D,\alpha):I$ for any general $\D\in\I2$ and any $\alpha\in I$ whose
image in $I/(\D)$ is general.
 \part 2 Then $\I2/I^2$ is cyclic, and  $\ell(\I i)=2$ for every even
$i>0$.
 \part 3 Then $\I i$ is a self-linked perfect ideal for every $i>0$.
 \part 4 Then the symbolic Rees algebra $\SR I$ is Noetherian; in fact,
$$\SR I=R[It,\D t^2].$$
 \pf
 Set $A:=R/(\D)$ and let `$\,\?{\phantom{s}}$' indicate the image in
$A$.  First, we show that, for every $i\ge2$,
	$$\I i = \alpha^{i-1}I + \bigl((\D)\cap\I i\bigr).\dno3.1$$
 That equation is equivalent to this one: $\?{\I i}={\?\alpha}^{i-1}\?I$.
The latter may be checked at each associated prime $p$ of the
$R$-module $\?{\I i}/{\?\alpha}^{i-1}\? I$, which is a submodule of
$A/{\?\alpha}^{i-1}\?I$.  So $p$ must be an associated prime of
$A/{\?\alpha}^{i-1}\?I$.  So $(A/{\?\alpha}^{i-1}\?I)_p$ has depth 0.
Hence $({\?\alpha}^{i-1}\?I)_p$ has depth 1.  Hence $\?I_p$ has depth
1.  Now, $A/\?I=R/I$, and so there is a short exact sequence,
$$0\to\?I_p\to A_p\to(R/I)_p\to0.$$
 Therefore, either $\depth A_p\le1$ or $\depth(R/I)_p=0$.  In the
former case, $\depth R_p\le2$.  Hence, again $\depth(R/I)_p=0$.  So, in
any case, either $I\not\subset p$, or $p$ is an associated prime of the
ideal $I$.  Hence $(\I i)_p=(I^i)_p$.  On the other hand, Parts (2) and
(4) of Proposition~\Cs1) yield ${\?I}^i= {\?\alpha}^{i-1}\?I$.  Hence
$(\?{\I i})_p= ({\?\alpha}^{i-1}\?I)_p$, and the proof of \Cs3.1) is
complete.

Next, we show that, for every $i\ge2$,
$$(\D)\cap\I i = \D\I{i-2}.\dno3.2$$
 First, consider the case where $I$ is generated by a regular sequence.
Suppose that $\D\in mI^2$, where $m$ denotes the maximal ideal of $R$.
Then
$$I^2=I^2\cap(\D,\alpha)=I^2\cap(mI^2,\alpha)
=mI^2+(\alpha)\cap I^2.$$
 So Nakayama's lemma implies that $I^2=(\alpha)\cap
I^2\subset(\alpha)$.  However, this inclusion is impossible because
$I^2$ has height 2.  Thus $\D\in I^2\setminus m I^2$.  Hence, in the
associated graded ring $\gr_IR$, the leading form of $\D$ is a regular
element, because $\gr_IR$ is a polynomial ring.  Consequently,
$(\D)\cap I^i = \D I^{i-2}$.  Thus \Cs 3.2) holds for complete
intersections.

To prove \Cs 3.2) in general, it suffices to show that $(\D)\cap\I
i$ lies in $\D\I{i-2}$.  So let $p$ be an associated prime of the ideal
$\D\I{i-2}$.  We may assume that $I\subset p$,  since the assertion is
trivial otherwise.  But, then
$$\depth(\I{i-2})_p=\depth(\D\I{i-2})_p=1.$$
 Hence $p$ is an associated prime of the ideal $\I{i-2}$ (and so
$i-2\ge1$).  Therefore, by the definition of symbolic powers, $p$ is
contained in the union of the associated primes of the ideal $I$; so
$p$ is contained in one of them.  Hence, again, $I_p$ is generated by a
regular sequence, and so the assertion follows from the discussion in
the preceding paragraph.

Together, \Cs 3.1) and \Cs 3.2) yield the following equation: for
$i\ge2$, $$\I i = \alpha^{i-1}I + \D\I{i-2}.\dno3.3$$
 In turn, that equation yields, via induction on $j$, this one: for
$j\ge1$,
	$$\I{2j}=\alpha(\D,\alpha^2)^{j-1}I+(\D^j).\dno3.4$$
 On the other hand, \Cs 3.3) immediately yields Assertion (4).

To prove Assertions (1) and (2), note that $\I2=\alpha I+(\D)$ by \Cs
3.3).  So $\I2/I^2$ is cyclic.  Furthermore,
	$$(\I2)^2=\alpha^2I^2+\D\alpha I+(\D^2)
		\subset(\D,\alpha^2)\I2\subset(\I2)^2.$$
 Hence $(\I2)^2=(\D,\alpha^2)\I2$.  Since $(\D,\alpha^2)$ lies in $\I2$,
it is therefore a reduction of $\I2$.  Hence $\ell(\I2)=2$.  Therefore,
Lemma~(3.2) implies Assertion (1).  Moreover, similarly, \Cs 3.4)
implies that $(\D^j,\alpha^{2j})$ is a reduction of $\I{2j}$; hence,
$\ell(\I{2j})=2$.  Thus Assertion (2) holds.

Together, \Cs 3.2) and \Cs 3.3) yield the following short exact sequence
for  $i\ge2$:
$$0\to\I{i-2}@>\D>>\I i\to{\?\alpha}^{i-1}\?I\to0.\dno3.5$$
 Since the $R$-module $I$ has projective dimension at most 1, so do
$\?I$ and therefore also ${\?\alpha}^{i-1}\?I$.  Hence it follow from
\Cs 3.5), via induction on $i$, that $\I i$ too has projective
dimension at most 1.  Therefore, $\I i$ is perfect, as asserted in (3).

To prove that $\I i$ is self-linked, we first treat the case where $i$
is odd, say, $i=2j-1$ with $j\ge1$.  In fact, we claim that
	$$\I{2j-1}=(\D^j,\alpha^{2j-1}):\I{2j-1}.\dno3.6$$
 Now, $I=(\D,\alpha):I$. So $I^2\subset(\D,\alpha)\subset I$.  Hence
$I$ and $(\D,\alpha)$ have the same radical.  Hence $I^{2j-1}$ and
$(\D^j,\alpha^{2j-1})$ have the same radical as $I$.  Therefore, it
suffice to establish \Cs 3.6) locally at each associated prime $p$ of
$I$.  However, by localizing at $p$, we may assume that $\I{2j-1}=
I^{2j-1}$, that $I=(\alpha,\beta)$, and that $\D=
r\alpha^2+s\alpha\beta+t\beta^2$ for some $r,s,t$.  Here $t$ has to be
a unit; otherwise, $\I2\neq\alpha I+(\D)$, contradicting \Cs 3.3).  So
we may suppose that $\D= r\alpha^2+s\alpha\beta+\beta^2$.

We thus need to show that
	$$(\D^j,\alpha^{2j-1}):(\alpha,\beta)^{2j-1}
		=(\alpha,\beta)^{2j-1}\dno3.7$$
 where
 $\D^j=\sum_{k=0}^{2j-2}b_k\alpha^{2j-1-k}\beta^k+\beta\beta^{2j-1}$
with $b_k\in(\alpha,\beta)$.  Now, $(\alpha,\beta)^{2j-1}$ is generated
by the signed maximal minors $d_1,\dots,d_{2j}$ of the following $2j$ by
$2j-1$ matrix:
 $$\pmatrix -\beta\cr
   \alpha&-\beta\cr
 &\alpha&\ddots\cr
 &&\ddots&-\beta\cr
 &&&\alpha\cr\endpmatrix.$$
 Adjoin the two columns of coefficients that arise when $\D^i$ and
$\alpha^{2j-1}$ are written in terms of the minors $d_1,\dots,d_{2j}$,
obtaining the following $2j$ by
$2j+1$ matrix:
$$\pmatrix -\beta&&&&b_0&1\cr
   \alpha&-\beta&&&b_1&0\cr
 &\alpha&\ddots&&\vdots&\vdots\cr
 &&\ddots&-\beta&b_{2j-2}&0\cr
 &&&\alpha&\beta&0\cr\endpmatrix.$$
 Its maximal minors generate the link
$(\D^j,\alpha^{2j-1}):(\alpha,\beta)^{2j-1}$;
 see \cite{AN, \p.316}.  The preceding matrix may be reduced to the
following one:
$$\pmatrix \alpha&-\beta&&&b_1\cr
 &\alpha&\ddots&&\vdots\cr
 &&\ddots&-\beta&b_{2j-2}\cr
 &&&\alpha&\beta\cr\endpmatrix.$$
 Denote the minor obtained by deleting the  $i$th column by $D_i$.  Then
$$(D_l,\dots,D_{2j})=(\alpha^{l-1}\beta^{2j-l},\dots,\alpha^{2j-1})$$
 for $1\le l\le 2j$, as can be proved easily via descending induction
on $l$ because $b_k\in(\alpha,\beta)$.  In particular,
$(D_1,\dots,D_{2j})=(\alpha,\beta)^{2j-1}$.  Thus the asserted equation
\Cs 3.7) holds.

If $i$ is even, say $i=2j$ with $j\ge1$, then we claim that
	$$\I{2j}=(\D^j,\alpha^{2j+1}):\I{2j}.$$
 Set $C:=R/(\D^j)$; let `$\,\wt{\phantom{s}}$' indicate the image in
$C$; and let $L$ denote the total ring of quotients of $C$.  In these
terms, our claim asserts that $\wt{\I{2j}}$ is equal to
$({\wt\alpha}^{2j+1}):_C\wt{\I{2j}}$.  However,
$$\eqalign{\wt{\I{2j}}&=\wt\alpha\wt{\I{2j-1}}\cr
 &=\wt\alpha\bigl(({\wt\alpha}^{2j-1}):_C\wt{\I{2j-1}}\bigr)\cr
 &=\wt\alpha\bigl(({\wt\alpha}^{2j-1}):_L\wt{\I{2j-1}}\bigr)\cr
 &=({\wt\alpha}^{2j+1}):_L(\wt\alpha\wt{\I{2j-1}})\cr
 &=({\wt\alpha}^{2j+1}):_C(\wt\alpha\wt{\I{2j-1}})\cr
 &=({\wt\alpha}^{2j+1}):_C\wt{\I{2j}}.}$$
 Indeed, the first and the last equations hold by \Cs 3.4); the second
equation holds by the case of an odd $i$; the third and fifth equations
hold because $\wt\alpha$ lies in $\wt I$ and is a nonzerodivisor; and
the fourth equation holds because $\wt\alpha$ is a nonzerodivisor.
Thus the assertion about self-linkage holds.  So Assertion (3) holds,
and the proof of \Cs3) is complete.

\thm4 Let $R$ be a Cohen--Macaulay local ring with infinite residue
field, and $I$ a perfect ideal of grade $2$ that is a complete
intersection at each associated prime.  Then the following conditions
are equivalent:\smallskip
 \item i $I=(\D,\alpha):I$ for some $\D\in\I2$ and some $\alpha\in I$;
 \item i$'$ $I=(\D,\alpha):I$ for any general $\D\in\I2$ and any
$\alpha\in I$ whose image in $I/(\D)$ is general;
 \item ii  $\I2/I^2$ is cyclic, and  $\ell(\I 2)=2$;
 \item iii $\SR I=R[It,\D t^2]$ for some $\D\in\I2$, and $\I i$ is a
perfect ideal for every $i>0$.
 \item iv $\SR I=R[It,\D t^2]$ for some $\D\in\I2$, and $\I i$ is a
Cohen--Macaulay ideal for infinitely many $i>0$;
 \smallskip\noindent
 Moreover, if one of the preceding conditions obtains, then $\I i$ is
self-linked for every $i>0$, and $\SR I$ is a Cohen--Macaulay ring.
 \pf
 First, (ii) implies (i$'$) by Lemma~\Cs 2), and (i$'$) implies (i)
trivially.  Second, (i) implies (ii) and (iii), and (i) implies that $\I
i$ is self-linked, all by Proposition~\Cs3).  (Those implications do not
require the Cohen--Macaulay hypothesis, but the remaining implications
do.)  Third, clearly, (iii) implies (iv).  Fourth, assume (iv).  Then
$\SR I=R[It,\I2 t^2]$; hence, $\I{2i}=(\I2)^i$ for every $i>0$.
Therefore, infinitely many powers of $\I2$ are Cohen--Macaulay ideals of
height 2.  Hence, $\ell(\I2)=2$; see \cite{Brod, Thm. 2, p. 36}.  On the
other hand, $\I2= I^2+(\D)$; so $\I2/I^2$ is cyclic.  Thus (ii) holds,
and so (i) to (iv) are equivalent.  Finally, (i) implies that $\SR I$ is
Cohen--Macaulay by the next proposition.

\prp5 Let $R$ be a Noetherian local ring, $I$ a perfect ideal of grade
$2$.  Assume that $I$ is a complete intersection at each associated
prime, and that $I=(\D,\alpha):I$ for some $\D\in\I2$ and some
$\alpha\in I$.  If $R$ is Cohen--Macaulay or Gorenstein, then so is
$\SR I$.
 \pf
 Set $A:=R/(\D)$ and let `$\,\?{\phantom{s}}$' indicate the image in
$A$.  Form the graded ring $G:=\bigoplus_{i\ge0}\I i/\I{i+1}$
associated to the filtration $\{\,\I i\mid i\ge0\,\}$ of $R$.  Denote
the leading form of $\D$ in $G$ by $\D'$.  Then $\deg\D'\ge2$ since
$\D\in\I2$.  On the other hand, for every $i\ge1$,
	$$\D\bigl(\I{i+1}:(\D)\bigr)=(\D)\cap\I{i+1}=\D\I{i-1};$$
  indeed, the first equation is obvious, and the second is \Cs3.2).
Hence
	$$\I{i+1}:(\D)=\I{i-1}$$
 because $\D$ is $R$-regular.  Therefore, $\deg\D'=2$  and $\D'$ is
$G$-regular.  Furthermore, again by \Cs3.2), for every $i\ge2$,
	$$\I{i+1} +\D\I{i-2} = \I{i+1} + \bigl((\D)\cap\I i\bigr).$$
  Therefore, $G/(\D')=\bigoplus_{i\ge0}\?{\I i}\bigm/\?{\I{i+1}}$.  Finally,
(3.3.3) gives $\?{\I i}=\?{I^i}$ for every $i\ge2$; in fact, then,
obviously, this equation holds for every $i\ge0$.  Therefore,
$$\textstyle G/(\D')=\bigoplus_{i\ge0}\?{\I i}\Bigm/\?{\I{i+1}}
=\bigoplus_{i\ge0}\?{I^i}\Bigm/\?{I^{i+1}}
=\gr_{\?I}A.$$
 Thus there is a short exact sequence of graded $G$-modules,
	$$0@>>>G(-2)@>\D'>>G@>>>\gr_{\?I}A@>>>0.\dno5.1$$

Let $\?\alpha'$ denote the leading form of $\?\alpha$ in $\gr_{\?I}A$.
Then $\?\alpha'$ is a regular element on $\gr_{\?I}A$ of degree one;
indeed, $\?{I^{i+1}}:(\?\alpha)=\?{I^i}$ because $\?\alpha$ is an
$A$-regular element contained in $\?I$ and because
${\?I}^{i+1}=\?\alpha{\?I}^i$ for every $i\ge1$ by \Cs3.3).  Moreover,
$\?\alpha{\?I}^{i-1}$ is equal to $\?\alpha\cap{\?I}^i$ because the
latter is equal to $\?\alpha\bigl(\?{I^{i}}:(\?\alpha)\bigr)$.
Therefore, $\gr_{\?I}A/(\?\alpha')$ can be identified with the
associated graded ring of the ideal $\wt I:=I/(\D,\alpha)$ in
$R/(\D,\alpha)$.  The latter graded ring is equal to the trivial
extension $S\ltimes \wt I(-1)$ of $S:=R/I$ because ${\wt I}\/^2=0$.
Thus there is an exact sequence of graded $G$-modules,
	$$0@>>>\gr_{\?I}A(-1)@>\?\alpha'>>\gr_{\?I}A
	@>>>S\ltimes \wt I(-1)@>>>0.\dno5.2$$
 Sequences \Cs5.1) and \Cs5.2) say, in other words, that the leading
forms $\D'$ and $\alpha'$ of $\D$ and $\alpha$ in $G$ have degrees two
and one, that they form a regular sequence, and that $G/(\D')$ is equal
to $\gr_{\?I}A$ and $G/(\D',\alpha')$ is equal to $S\ltimes \wt I(-1)$.

Set $d:=\dim R$.  Let $\m $ denote the maximal ideal of $R$, and
$\bold M$ the irrelevant maximal ideal of $G$.  Given a finitely
generated graded $G$-module of dimension $s$, let $H_\bold M^j(E)$
denote the $j$th local cohomology module of $E$ with supports in $\bold
M$, and set
	$$a(E):=\sup\{\,i\mid[H_\bold M^s(E)]_i\not=0\}.$$

Assume now that $R$ is a Cohen--Macaulay ring.  Then $S$ is a
Cohen--Macaulay ring, and $\wt I$ a maximal Cohen--Macaulay $S$-module
because $I$ is a perfect $R$-ideal of grade 2; hence, $S\ltimes
\wt I(-1)$ is a Cohen--Macaulay ring.  Therefore, $\gr_{\?I}A$ and $G$ are
Cohen--Macaulay too, because $\D'$ and $\alpha'$ form a regular
sequence.  Furthermore,
$$\eqalign{
 H_\bold M^{d-2}(S\ltimes \wt I(-1))&=H_{\m G}^{d-2}(S\ltimes
       \wt I(-1))=H_\m ^{d-2}(S\ltimes \wt I(-1))\cr
	&=H_\m ^{d-2}(S)\oplus H_\m ^{d-2}(\wt I)(-1).\cr}$$
 Hence $a(S\ltimes \wt I(-1))=1$.  Hence \Cs5.2) implies
$a(\gr_{\?I}A)=0$ because $\gr_{\?I}A$ is Cohen--Macaulay.  Hence
\Cs5.1) implies $a(G)=-2$ because $G$ is Cohen--Macaulay.  Now, by
\cite{G-N, (1.2), \p.74}, since $a(G)<0$ and since $G$ is
Cohen--Macaulay, $\SR I$ is Cohen--Macaulay.

Suppose $R$ is Gorenstein.  Then the dualizing module $\omega_S$ of $S$
is isomorphic to $\Ext_R^2(S,R)$.  However, $\D,\alpha$ form a regular
sequence, so this Ext is equal to
$\Hom_R(S,R/\bigl(\D,\alpha)\bigr)$ or, what is the same, to
$$\bigl((\D,\alpha):I\bigr)/(\D,\alpha).$$
   This module is equal to $\wt I$.  Thus $\wt I\cong\omega_S$.  So
$S\ltimes\wt I\cong S\ltimes\omega_S$.  The latter is a Gorenstein ring
(see \cite{Reiten, (7), \p.419} for instance).  Hence $G$ is Gorenstein
too, because $\D'$ and $\alpha'$ form a regular sequence.  Finally, by
\cite{G-N, (1.4), \p.75}, since $a(G)=-2$ and since $G$ is Gorenstein,
$\SR I$ is Gorenstein.

\prp6
 Let $R$ be a Noetherian local ring, and $B$ an $R$-algebra.  Assume
that, as an $R$-module, $B$ is finitely generated and perfect of grade
$1$.  Set $F_i:=\Fitt_i^R(B)$.  Set $I:=F_1$, and assume  $I\not= R$
and  $\grade I\ge2$.
Set $A:=R/F_0$ and let `$\,\?{\phantom{s}}$' indicate the image in $A$.
Given generators $1$, $u_2$, \dots, $u_n$ of the $R$-module $B$, let
$\varphi$ be an $n$ by $n$ matrix with entries in $R$ whose transpose
$\varphi^*$ presents $B$ via those generators.  Set $\D:=\det\varphi$.
Let $\phi$ be the $n$ by $n-1$ matrix consisting of the last $n-1$
columns of $\varphi$, and set $I':=\BI_{n-1}(\phi)$.
 \part 1 Then $I=I'$; moreover, $I$ is a perfect $R$-ideal of grade $2$,
and is presented by $\phi$.  Furthermore, $\Delta$ is regular.  In
addition, if $K$ denotes the total ring of quotients of $A$, then
$$A\subset B\subset K,\ \?I=A:_KB, \and B=A:_K\?I.$$
 \part 2 Assume that $\grade F_2\ge 3$.  Then $I$ is a complete
intersection at each associated prime, and $\D\in I^{(2)}$.
 \part 3 Assume that $B$ is Gorenstein over $R$, and modify
$\varphi$ to make it symmetric (see {Theorem~(2.3)}).  Let $\phi'$ be
the $n-1$ by $n-1$ matrix obtained by deleting the first row and column
of $\varphi$, and set $\alpha:=\det\phi'$ and set
$F_2':=\BI_{n-2}(\phi')$.  Then $\D,\alpha$ form an $R$-regular
sequence, with respect to which $I$ is self-linked.  Moreover,
${\?I}^{i+1}= \?\alpha{\?I}^i$ for every $i\ge1$, and
$B={\?\alpha\,}^{-1}\?I$.  In addition,
	$$F_2=\Fitt_2^R(I)= I^2:(\D)=F_2'.$$
 \pf
 Consider (1).  First of all, the proof of (2.8) yields, mutatis
mutandis, that $\sqrt I=\sqrt{I'}$.  Hence $\grade I'\ge2$.  Therefore,
by the Hilbert--Burch theorem, $I'$ is a perfect ideal of grade 2, and
is presented by $\phi$.  Second, $\varphi^*$ presents $B$, so
$B=\Cok(\?\varphi^*)$.  Moreover, $B$ is perfect of grade $1$; hence
$\varphi^*$ is injective and so $\Delta$ is regular.  Finally, consider
the proof of \Cs1)(1).  The latter half applies, therefore, to the
present $\varphi$, and yields $B=A:_K\?{I'}$ and $\?{I'}=A:_KB$.
Hence Proposition~\Cs1)(2) yields $I=I'$.  Thus (1) is proved.

Assertion (2) follows from Parts ($2'$) and (3) of Proposition~\Cs1),
which applies to the present $I$ and $\D$ because of (1).

Consider (3).  The first assertions follow, thanks to (1), from
Lemma~(2.5)(3) and the equivalence of (vi), (iv), and (ii) of
Proposition~(3.1)(4).  Now, $\phi$ presents $I$ by (1); hence,
$$F_2'\subset \Fitt_2^R(I)\subset F_2.$$
 On the other hand, we're about to show that
$$F_2\subset F_1^2:(\D)\subset F_2'.\dno6.1$$
 Then the displayed equations asserted in (3) will follow immediately.

As in the proof of Proposition~\Cs1)(2'), the Laplace identity (2.4)
yields
	$$F_2(\D)\subset F_1^2.$$
 Furthermore, as in the proof of Lemma~(2.5)(3), the Laplace identity
yields
	$$m_{1,k}^{j,1}\D=m_j^1m_k^1-\alpha m_k^j$$
 because $\varphi$ is symmetric.  Hence,
	$${I'}^2 \subset F_2'(\D)+(\alpha)F_1;\dno6.2$$
 Since $F_1=I'$ by (1), therefore
   $$F_2\subset F_1^2:(\D)={I'}^2:(\D)\subset (F_2'\D,\alpha):(\D).$$
 Now, it follows easily from the definitions that
$$(F_2'\D,\alpha):(\D)\subset{\biglp(\alpha):(\D)\bigr)}+F_2'.$$
 However, since $\alpha,\D$ form a regular sequence,
$$\bigl((\alpha):(\D)\bigr)=(\alpha).$$
 Together, the last three displays yield that
$$F_2\subset F_1^2:(\D)\subset (\alpha)+F_2' = F_2'.$$
 Thus  \Cs6.1) holds, and so (3) is proved.

\thm7 Let $R$ be a Noetherian local ring.  Given an element $\D\in R$,
set $A:=R/(\D)$ and let `$\,\?{\phantom{s}}$' indicate the image in
$A$, and let $K$ denote the total ring of quotients of $A$.
Given a finitely generated $R$-module $B$, set $F_i:=\Fitt_i^R(B)$.
Then there is a bijective correspondence between the following
two sets:
  $$\alignat2 \{B&\mid  B\text{ a perfect $R$-algebra of grade 1 with
	$\grade	F_1=2$ and $\grade F_2\ge3$}&&\},\\
	\{(\D)\subset I&\mid I\text{ a perfect $R$-ideal
		of grade $2$ that is a complete intersection}\\
	&& \llap{at each associated prime,
		and $\D\in\I2$ an $R$-regular element}&\}.\endalignat$$
 The correspondence associates to $B$ the pair $(\D):=F_0$ and
$I:=F_1$, and conversely to $(\D)\subset I$ the $R$-module
$B:=A:_K\?I$; moreover, $A\subset B\subset K$ and $\?I=A:_KB$.
 Finally, $B$ is Gorenstein if and only if $I=(\D,\alpha):I$ for some
$\alpha\in I$; if so, then $B={\?\alpha}^{-1}\?I$.
 \pf
 The assertions follow immediately from Propositions~(3.1) and (3.6)(3).

\prp8
 Let $R$ be a Noetherian local ring, and $B$ a Gorenstein $R$-algebra.
Assume that, as an $R$-module, $B$ is finitely generated and perfect of
grade $1$.  Set $F_i:=\Fitt_i^R(B)$.  Then either $\height F_2\le3$ or
$F_2=R$.

In addition, assume $\grade F_1\ge2$ and $\grade F_2\ge 3$.  Then $F_2$
is a perfect $R$-ideal of grade $3$, or it's the unit ideal.

Finally, let $\Delta$ be a generator of $F_0$, and set $I:=F_1$.  Then
$\D\in I^{(2)}$, and $\Delta$ defines an isomorphism,
	 $$R/F_2 \risom I^{(2)}/I^2.$$
 \pf
 Proposition~(2.9)(1) yields $\height F_2\le3$.  Proposition~\Cs6)(3)
yields a symmetric matrix $\phi'$ whose ideal of submaximal minors is
equal to $F_2$.  Hence $F_2$ is a perfect $R$-ideal of grade $3$ since
$\grade F_2\ge 3$; see \cite{Kutz, Thm.~1, \p.116} for instance.
Finally, by Proposition~(3.6)(2) and (3), the relation $\D\in I^{(2)}$
holds and the hypotheses of Proposition~(3.3) obtain.  So \Cs 3.3)
yields $I^{(2)}=I^2+(\D)$.  Hence, multiplication by $\Delta$ induces a
surjection $R\onto I^{(2)}/I^2$.  Obviously, its kernel is $I^2:(\D)$.
However, the latter is equal to $F_2$ by Proposition~\Cs 6)(3).

\cor9
 Let $f\:X\to Y$ be a finite map of locally Noetherian schemes.  Assume
that $f$ is locally of flat dimension $1$ and Gorenstein.  Then each
component of $N_3$, the scheme of target triple points, has codimension
at most $3$.

Assume in addition that each component of $N_3$ has codimension $3$,
that $f$ is birational onto its image, and that $Y$ satisfies Serre's
condition \Ser3.  Let ${\Cal I}_{N_2}$ denote the $\O_Y$-ideal of the
scheme of target double-points, and ${\Cal I}_{N_2}^{(2)}$ its second
symbolic power.  Then $\O_{N_3}$ is a perfect $\O_Y$-module, and ${\Cal
I}_{N_2}^{(2)}/{\Cal I}_{N_2}^2$ is an invertible
$\O_{N_3}$-module.
 \pf
 The assertions are local on $Y$.  Furthermore, $\O_Y$ is perfect of
grade 1 by virtue of \cite{KLU, Lemma~(2.3)}, and if $f$ is birational
onto its image and $Y$ satisfies \Ser2, then ${\Cal I}_{N_2}$ has grade
at least 2 by Lemma~(2.11).  So the assertions follow from
Proposition~(3.8).

\newsect Normality

In this section, we study the normality of the algebra $B$ treated in
Proposition~(3.1) and Theorem~(3.7).  In other words, we investigate
when $B$ is equal to the integral closure $A'$ of $A$ in its total ring
of quotients $K$.

Given a Noetherian local ring $R$, let $e(R)$ denote the multiplicity
of $R$, and call $R$ a hypersurface ring if $R\cong S/(x)$ where $S$ is
a regular local ring and $x$ is an $S$-regular element.  Given an ideal
$I$ and an element $x$, set
	$$O(I):=\sup\{i|I\subset m^i\} \and o(x):=O((x))$$
 where $m$ is the maximal ideal, and call these numbers the {\it
orders\/} of $I$ and $x$.

\prp1
 Let $R$ be a Noetherian local ring satisfying (R$_2$) and (S$_3$).  Let
$B$ be an $R$-algebra that, as an $R$-module, is finitely generated
and perfect of grade $1$, and assume that one of the equivalent
conditions of Lemma~(2.11) is satisfied.  Set $F_i:=\Fitt_i^R(B)$ and
$I:=F_1$.  Set $A:=R/F_0$ and let $\Delta$ be a generator of $F_0$.
Consider the following conditions: \smallskip
 \item i $e(A_p)=2$ for every associated prime $p$ of $I$, and
$\Delta\in\I2$;
 \item ii $e(A_p)\le3$ for every associated prime $p$ of $I$, and
$\Delta\in\I2$;
 \item iii $(R/I)_p$ is a hypersurface ring for every associated prime
$p$ of $I$;
 \item iv $I_p$ is a complete intersection for every associated prime
$p$ of $I$;
 \item v $\Delta\in\I2$.
 \smallskip\noindent
 Then (i)$\impls$(ii)$\impls$(iii)$\impls$(iv)$\impls$(v), and all five
conditions are equivalent if $B$ is normal.
 \pf
 Localizing at an associated prime of $I$, we may assume that $R$ is a
regular local ring of dimension 2, that $I$ is primary for the maximal
ideal $\m$, and that $\I2=I^2$.  First notice that $e(A_p)\ge2$, because
otherwise $A$ would be regular, and hence, since the extension $A\subset
B$ is finite and birational, $A$ would be equal to $B$, contrary to the
assumption that $I$, or $F_1$, is $\m$-primary, so unequal to $R$.

Trivially, (i) implies (ii).  Assume (ii).  If $I\subset\m^2$, then
$\D\in I^2\subset\m^4$ and so $e(A)=o(\D)\ge4$.  It follows that (iii)
holds.  Clearly, (iii) implies (iv).  By Lemma~(3.6)(1), the hypotheses
of Proposition~(3.1) obtain, and its Part~(3) yields (iv)$\impls$(v).

Finally, suppose that $B$ is normal and that $\Delta$ is in $\I2$, so in
$I^2$.  We have to prove that $e(A)=2$.  Set $\?I:=I/(\Delta)$.  Then
$\?I$ is equal to the conductor in $A$ of $B$ by Proposition~(3.1)(1).
So, since $B$ is one-dimensional and normal, $\?I\cong B$.  Since
$\Delta\in\m I$, therefore $$I/\m I\cong\?I/\m\?I\cong B/\m B.$$
 Now, the length of $I/\m I$ is equal to the minimal number of
generators $\nu(I)$ of $I$, and the length of $B/\m B$, viewed as an
$A$-module, is equal to the multiplicity of $e(A)$.  Hence
	$$\nu(I)=e(A)=o(\D).$$
 On the other hand, the Hilbert--Burch theorem immediately yields the
(well-known) inequality, $O(I)\ge\nu(I)-1$.  Also, $o(\D)\ge 2O(I)$
since $\Delta\in I^2$.  So
	$$o(\D)\ge 2O(I)\ge2\nu(I)-2\ge\nu(I).$$
 Together, the two preceding displays give $\nu(I)=2\nu(I)-2$.  So
$\nu(I)=2$, and the proof is complete.

\eg2 Proposition~(4.1) shows that the five conditions of
Proposition~(3.1)(2) are not equivalent in general.  Indeed, let $R$ be
a power series ring in two variables over a field (of any
characteristic), and let $\D$ be any reduced power series of order
$o(\D)$ at least 3, and set $A:=R/(\D)$.  Let $B$ be the integral
closure of $A$ in its total ring of quotients.  Then, by Lemma (2.11)
and Proposition~(3.6), the hypotheses of Proposition~(3.1) obtain.
However, $\D\notin\I2$ by virtue of Proposition~(4.1), although $B$ is a
ring.  For instance, say $R:=k[[x,y]]$ and $\D:=y^3-x^4$; then
$I=(x^2,xy,y^2)$ because, if $B=k[[t]]$ say, then $A=k[[t^3,t^4]]$ and
$I/(\D)=(t^6,t^7,t^8)$.

\lem3 Let $k$ be an infinite field, $L$ a field containing $k$.  Let
$\psi$ be an $s$ by $r$ matrix with coefficients in $L$ and with rank
$s$.  Then a general $k$-linear combination of the columns of $\psi$ does
not satisfy any given nonzero polynomial $F$ in $L[X_1,\dots,X_s]$.
 \pf
 Consider the polynomial $G:=F\circ\psi$ in $L[Y_1,\dots,Y_r]$.
Obviously, $G\neq0$ because $\rk\psi=s$ and $F\neq0$.  Let $\{a_i\}$ be
a basis for $L/k$, and say $G=\sum a_iG_i$ where $G_i\in
k[Y_1,\dots,Y_r]$.  Then $G_i\neq0$ for some $i$ because $G\neq0$.
Hence, since $k$ is infinite, there exists a $y\in k^r$ such that
$G_i(y)\neq0$ for some $i$.  On the other hand, for any such $y$,
clearly $G(y)\neq0$.  Thus the assertion holds.

\thm4 Let $k$ be an infinite perfect field, and $R$ a Noetherian local
$k$-algebra satisfying (R$_2$) and (S$_3$).  Let $I$ be a perfect ideal
of grade $2$ that is a complete intersection at each associated prime.
Fix generators $h_1,\dots,h_l$ of $\I2$, and generators
$x_1,\dots,x_n$ of the maximal ideal of $R$.  Let $f_1,\dots,f_r$ be
the sequence $h_1,\dots,h_l$ if $\char k=0$, and the sequence
$x_0h_1,\dots,x_nh_l$ with $x_0:=1$ if $\char k$ is arbitrary.  Let
$\D$ be a general $k$-linear combination of $f_1,\dots,f_r$, and set
$A:=R/(\D)$.  Let $K$ be the total ring of quotients of $A$, and set
$B:=A:_K\?I$.  Then $B$ is a ring, and the following conditions are
equivalent:
 \smallskip
 \item i $e(A_p)=2$ for every associated prime $p$ of $I$;
 \item ii $e(A_p)\le3$ for every associated prime $p$ of $I$;
 \item iii $(R/I)_p$ is a hypersurface ring for every associated prime
$p$ of $I$;
 \item iv $B$ is normal.
 \smallskip\noindent
 Furthermore, if any of these conditions hold, then the extension $B/A$
is unramified in codimension 1.
 \pf
  First of all, $B$ is a ring by Proposition~(3.1)(2), and $B$ is a
perfect $R$-module of grade 1 by Proposition~(3.1)(1).  So we may apply
Proposition~(4.1).  Hence, (i)$\impls$(ii)$\impls$(iii); moreover,
(iv)$\impls$(i) since $\Delta\in\I2$.  Finally, assume (iii).  We have
to prove that $B$ is normal and $B/A$ is unramified in codimension 1.

By hypothesis, $R$ satisfies (S$_3$).  Since $B$ is a perfect
$R$-module of grade 1, it is locally of codimension $1$ by \cite{KLU,
(2.3) and (2.5)}.  Hence $B$ satisfies (S$_2$) as an $R$-module, and so
as a ring.  Therefore, it suffices to prove this: for every prime $p$
of $R$ with $p\ni\D$ and with $\dim R_p\le2$, the localization $B_p$ is
regular and is unramified over $A_p$.

If $p\not\supset I$, then $A_p$ is regular by a form of Bertini's
theorem since $R_p$ is regular and since $\D\in\I2$ is general.
Indeed, in characteristic zero, the required form is \cite{Flenner,
Thm.~4.6, \p.107}; in arbitrary characteristic, the required form
follows from the proof of \cite{Flenner, Thm.~4.1, \p.106} and
\cite{Flenner, Thm.~1.7, \p.101}.  Thus we may assume that $p\supset
I$.  Then $p$ is one of the finitely many associated primes of $I$.

Localize at $p$.  Then $R$ is a regular local ring of dimension 2, and
$I\not\subset\m^2$ where $\m$ is the maximal ideal because $R/I$ is a
hypersurface ring.  Completing $R$, we may even suppose that
$R=L[[x,y]]$ where $L$ is an extension field of $k$ and $x,y$ are
indeterminates with $x\in I$.  Then $I=(x,y^s)$ for some $s\ge1$.  Set
$\?I:=I/(\D)$.  Let $\D_x$ stand for the partial derivative of $\D$
with respect to $x$, and $D$ for the image of $\D_x$ in $A$.  We are now
going to establish the following equation:
	$$D\?I=\?I^2.\dno4.1$$

Since both $f_1,\dots,f_r$ and $x^2,xy^s,y^{2s}$ generate $I^2$, there
exist a 3 by $n$  matrix $\varphi$  and a $n$ by 3 matrix $\psi$, both
with entries in $R$, so that
	$$\col x^2 xy^s y^{2s} = \varphi\col f_1 {\vdots} f_r
	\and \col f_1 {\vdots} f_r = \psi\col x^2 xy^s y^{2s} .$$
 Since $x^2,xy^s,y^{2s}$ form a minimal generating set, $\varphi\psi$
is congruent to the 3 by 3 identity matrix modulo $m$.  Hence $\psi$
has rank 3 modulo $m$; in other words,
	$$\rk_L\psi(0)=3$$
  where $\psi(0)$ is the matrix of constant terms.  By definition,
	$$\D:=a_1f_1+\cdots+a_rf_r$$
 where $(a_1,\dots,a_r)$ is a general point in $\bold{A}_k^n(k)$.  Set
	$$(u,v,w):=(a_1,\dots,a_r)\psi.$$
 Then the vector of constant terms $(u(0),v(0),w(0))$ is a general
$k$-linear combination of the rows of $\psi(0)$; hence, it does not
satisfy any given nonzero polynomial with coefficients in $L$ by
Lemma~\Cs3) applied to the transpose of $\psi(0)$.  On the other hand,
	$$\D = ux^2+vxy^s+wy^{2s}.$$

Set $\wt u:=2u+xu_x$ and $\wt v:=v+xv_x+y^sw_x$.  Then
	$$\D_x=\wt ux +\wt vy^s.$$
  Since $(u(0),v(0),w(0))$ satisfies no given nonzero polynomial,
$w(0)\neq0$ and
    $$({\wt u}^2w-\wt uv\wt v+u{\wt v}^2)(0)
	= 4u^2(0)w(0)-u(0)v^2(0)\neq0.$$
 So $w$ and ${\wt u}^2w-\wt uv\wt v+u{\wt v}^2$ are units in $R$.
Since $w$ is a unit,
	$$\align(\D,\D_xI)
   &=(ux^2+vxy^s+wy^{2s},(\wt ux+\wt vy^s)x,w(\wt ux+\wt vy^s)y^s)\\
&=(ux^2+vxy^s+wy^{2s},\wt ux^2+\wt vxy^s,w\wt uxy^s-\wt v(ux^2+vxy^s))\\
&=(ux^2+vxy^s+wy^{2s},\wt ux^2+\wt vxy^s,-\wt vux^2+(w\wt u-\wt vv)xy^s).\\
   \endalign$$
 However,
   $$\det\pmatrix \wt u& \wt v\\-\wt vu &(w\wt u-\wt vv)\\\endpmatrix
	={\wt u}^2w-\wt uv\wt v+u{\wt v}^2,$$
 and the later is a unit.  Hence
	$$(\D,\D_xI)=(\D,x^2,xy^s)=(\D,x^2,xy^s,y^{2s})=(\D,I^2).$$
 So \Cs4.1) holds.

Since $D\?I=\?I^2$ by \Cs4.1), Proposition~(3.1)(4) yields that
$\?I=DB$ and that $D$ is a nonzerodivisor.  The latter implies that the
extension $A/L[[y]]$ is generically unramified.  Hence $A$ is
generically reduced, so reduced since it has no embedded components.
Therefore, the integral closure $A'$ of $A$ in $K$ is a finitely
generated $A$-module.

Let $C$ denote the conductor.  Then $C$ contains $D$, and $C=DA'$ if
and only if $A'/L[[y]]$ is unramified; this is a matter of standard
theory, see \cite{Kunz, Cor.~G.12, \p.389} and \cite{Kunz, Cor.~G.14c),
\p.391} for example.  On the other hand, $B\subseteq A'$ since $B$ is a
finitely generated $A$-module; hence, $C\subseteq\?I$ because
$\?I=A:_KB$ by Proposition~(3.1)(1).  Therefore,
	$$C\subseteq \?I=DB\subseteq DA'\subseteq C.$$
 So $DB=DA'$; whence $B=A'$, and thus $B$ is normal.  Moreover,
$DA'=C$; whence $A'/L[[y]]$ is unramified, and so $A'/A$ is too.  The
proof is now complete.

\eg5 In Theorem~\Cs4), it is necessary to assume that $\D$ is general.
For example, let $k$ be a field of any characteristic, and
$R:=k[[x,y]]$ a power series ring.  Set $I:=(x,y)$ and $\D:=x^2-y^5$.
Then $A'$ is a power series ring in one variable, say $A'=k[[t]]$, and
$A$ and $B$ are subrings; namely, $A=k[[t^5,t^2]]$ and
$B=k[[t^3,t^2]]$.  So $\?I=(t^5,t^2)$ and $C=(t^5,t^4)$.  Moreover,
$D=t^5$ if $\char k\neq2$, and $D=0$ if $\char k=2$.

\newsect References

\references

AN
 M. Artin and M. Nagata
 \paper Residual intersections in Cohen--Macaulay rings
 \jmku 12 1972 307--23

Brod
 M. Brodmann
 \paper The asymptotic nature of analytic spread
 \mpcps 86 1979 35--39

Bruns
 W. Bruns
 \paper The Eisenbud--Evans generalized principal i\-de\-al theorem and
determinantal i\-de\-als
 \pams 83 1981 19--24

BE77a
 D. Buchsbaum and D. Eisenbud
 \paper What annihilates a module
 \ja 47 1977 231--43

BE77b
 D. Buchsbaum and D. Eisenbud
 \paper Algebraic structures for finite free resolutions, and some
structure theorems for ideals of codimension 3
 \ajm 99 1977 447--85

Cat
 F. Catanese
 \paper Commutative algebra methods and equations of regular surfaces
 \bucharest 68--111

dJvS
 T. de Jong and D. van Straten
 \paper Deformations of the normalization of hypersurfaces
 \ma 288 1990 527--47

EN
 J. Eagon and D. Northcott
 \paper Ideals defined by matrices and a certain complex associated to
them
 \prs A269 1962 188--204

Eisenbud
 D. Eisenbud
 \paper Homological algebra on a complete intersection, with an
application to group representations
 \tams 260 1980 35--64

Ei-Ma
 D. Eisenbud and B. Mazur
 \paper Symbolic squares, Fitting ideals, and evolutions
 \paperinfo Preprint\egroup

Flenner
 H. Flenner
 \paper Die S\"atze von Bertini f\"ur lokale Ringe
 \ma 229 1977 97--111

G-N
 S. Goto and K. Nishida
 \paper The Cohen--Macaulay and Gorenstein Rees algebras associated to
filtrations
 \mams 110 1994

G-N-S
 S. Goto, K. Nishida, and Y. Shimoda
 \paper The Gorensteinness of symbolic Rees algebras for space curves
 \jmsj 43 1991 465--81

G-N-W
 S. Goto, K. Nishida, and K. Watanabe
 \paper Non-Cohen--Macaulay symbolic blow-ups for space monomial curves
 and counterexamples to Cowsik's question
 \pams 120 1994 383--92

Grassi
 M. Grassi
 \paper Koszul modules and Gorenstein algebras
 \paperinfo Preprint, Universit\`a di Pisa, Jan.~1995\egroup

Groth
 A. Grothendieck
 \paper Cohomologie locale des faisceaux coh\'erents et th\'eor\`emes de
Lefschetz locaux et globaux
 \paperinfo (SGA II), Expos\'e XIV, Advanced Studies in Math, Vol. II,
North Holland Publishing Co., Amsterdam, 1968\egroup

HeK
J. Herzog and M. K\"uhl
 \paper Maximal Cohen--Macaulay modules over Gorenstein rings and
Bourbaki sequences
 \paperinfo in ``Commutative Algebra and Combinatorics,'' M. Nagata and
H. Matsumura (eds.), Advanced Studies in Pure Math {\bf 11},
North-Holland, Amsterdam (1987), 65--92\egroup

H-U
 J. Herzog and B. Ulrich
 \paper Self-linked curve singularities
 \nmj 120 1990 129--53

Jzfk
 T. Joz\'efiak
 \paper Ideals generated by minors of a symmetric matrix
 \cmh 53 1978 594--607

JP
 T. Joz\'efiak and P. Pragacz
 \paper Ideals generated by Pfaffians
 \ja 61 1979 189--98

KLU92
  S. Kleiman, J. Lipman, and B. Ulrich
 \paper The source double-point cycle of a finite map of codimension one
 \paperinfo in ``Complex Projective Varieties,'' G. Ellingsrud, C.
Peskine, G. Sacchiero, and S. A. Stromme (eds.)
 \lmslns 179 1992 199--212

KLU
 S. Kleiman, J. Lipman, and B. Ulrich
 \paper The multiple-point schemes of a finite curvilinear map of
codimension one\egroup

Kunz
 E. Kunz
 \book K\"ahler differentials
 \bookinfo  Advandced Lectures in Math.,
 Vieweg, 1986

Kutz
 R. Kutz
 \paper Cohen--Macaulay rings and ideal theory of invariants of
algebraic groups
 \tams 194 1974 115--29

Mats
 H. Matsumura
 \book Commutative ring theory
 \bookinfo Cambridge studies in advanced math {\bf 8}, 1986

MP
 D. Mond and R. Pellikaan
 \paper Fitting ideals and multiple points of analytic mappings
 \patzcuaro 107--161

Reiten
 I. Reiten
 \paper The converse of a theorem of Sharp on Gorenstein modules
 \pams 40 1972 417--20

Roberts85
 P. Roberts
 \paper A prime ideal in a polynomial ring whose symbolic blow-up is
not Noetherian
 \pams 94 1985 589--92

Roberts90
 P. Roberts
 \paper An infinitely generated symbolic blow-up in a power series ring
and a new counterexample to Hilbert's fourteenth problem
 \ja 132 1990 461--73

Shamash
 J. Shamash
 \paper The Poincar\'e series of a local ring
 \ja 12 1969 453--70

VV
 P. Valabrega and G. Valla
 \paper Form rings and regular sequences
 \nmj 72 1978 93--101

Valla84
 G. Valla
 \paper On set-theoretic complete intersections
 \arcireale 85--101

\endreferences

\enddocument